\documentclass[ejs]{imsart}
\doi{10.1214/154957804100000000}
\pubyear{0000}
\volume{0}
\firstpage{0}
\lastpage{0}\arxiv{1308.0417}

\usepackage{amsthm,amsmath,verbatim}
\RequirePackage[numbers]{natbib}
\RequirePackage[colorlinks,citecolor=blue,urlcolor=blue]{hyperref}

\usepackage{amssymb}
\usepackage{dsfont}
\usepackage{color}

\startlocaldefs

\newcommand\indicator{\mathds{1}}
\newcommand\prob{\mathds{P}}
\newcommand{\md}{\mathrm{d}}
\newcommand{\R}{\mathbb{R}}
\newcommand{\E}{\mathbb{E}}
\newcommand{\eps}{\varepsilon}

\newtheorem{theorem}{Theorem}[section]
\newtheorem{corollary}{Corollary}[section]
\newtheorem{lemma}{Lemma}[section]
\newtheorem{remark}{Remark}[section]
\newcommand{\mono}{(A1)}
\newcommand{\regul}{(A4)}
\newcommand{\positive}{(A5)}
\newcommand{\Bn}{(A2)}
\newcommand{\Bnn}{(A3)}
\newcommand{\tqed}{\hfill{$\Box$}}
\renewcommand{\hat}{\widehat}

\endlocaldefs

\begin{document}
\begin{frontmatter}

\title{A Kiefer-Wolfowitz type of result in a general setting, with an application to smooth monotone estimation}
\runtitle{A Kiefer-Wolfowitz type of result}

\author{\fnms{C\'{e}cile} \snm{Durot}\corref{}\ead[label=e1]{cecile.durot@gmail.com}}
\address{Universit\'e Paris Ouest Nanterre La D\'efense, Nanterre, France\\ \printead{e1}}
\and
\author{\fnms{Hendrik P.} \snm{Lopuha\"a}\ead[label=e2]{H.P.Lopuhaa@tudelft.nl}}
\address{Delft University of Technology, Delft, The Netherlands\\ \printead{e2}}
\runauthor{C. Durot and H.P. Lopuha\"a}



\begin{abstract}
We consider Grenander type estimators for monotone functions $f$ in a very general setting, which includes
estimation of monotone regression curves, monotone densities, and monotone failure rates.
These estimators are defined as the left-hand slope of the least concave majorant~$\widehat{F}_n$ of
a naive estimator $F_n$ of the integrated curve $F$ corresponding to~$f$.
We prove that the supremum distance between $\widehat{F}_n$ and $F_n$ is of the order
$O_p(n^{-1}\log n)^{2/(4-\tau)}$,
for some $\tau\in[0,4)$ that characterizes the tail probabilities of an approximating process for $F_n$.
In typical examples, the approximating process is Gaussian and $\tau=1$, in which case the convergence rate is $n^{-2/3}(\log n)^{2/3}$
is in the same spirit as the one obtained by Kiefer and Wolfowitz~\cite{kiefer-wolfowitz1976} for the special
case of estimating a decreasing density.
We also obtain a similar result for the primitive of $F_n$, in which case $\tau=2$,  leading to
a faster rate $n^{-1}\log n$, also found by Wang and Woodfroofe~\cite{wang2007kiefer}.
As an application in our general setup, we show that a smoothed Grenander type estimator
and its derivative are asymptotically equivalent to the ordinary kernel estimator
and its derivative in first order.
\end{abstract}





\end{frontmatter}

\section{Introduction}
\label{sec:intro}
Grenander~\cite{grenander1956theory} proved that the maximum likelihood estimator of a distribution $F$ that is concave on its support,
is the least concave majorant~$\widehat{F}_{n}$ of the empirical distribution function~$F_{n}$ of the $n$ independent observations.
In the case where~$F$ is absolutely continuous with probability density function $f$, the concavity assumption on $F$ simply means that
$f$ is non-increasing on its support, and the so-called Grenander estimator of $f$ is the left-hand slope of $\widehat F_{n}$.
Kiefer and Wolfowitz~\cite{kiefer-wolfowitz1976} showed that~$\widehat F_{n}$ and $F_{n}$ are close for large $n$ and as a consequence,
that~$\widehat F_{n}$ enjoys similar optimality properties as $F_{n}$, with the advantage of taking care of the shape constraint of being concave.
Roughly speaking, Kiefer and Wolfowitz~\cite{kiefer-wolfowitz1976} prove in their Theorem 1 that,
if $f$ is bounded away from zero with a continuous first derivative $f'$ that is bounded and bounded away from zero,
then, with probability one, the supremum distance between~$\widehat{F}_{n}$ and $F_{n}$ is of the order $n^{-2/3}\log n$.
Their main motivation was to prove the asymptotic minimax character  of $\widehat F_{n}$.
Their result easily extends to the case of an increasing density function,
replacing the least concave majorant with the greatest convex minorant.

In the setting of estimating an increasing failure rate,
Wang~\cite{wang1986asymptotically} proves that under appropriate assumptions,
the supremum distance between the empirical cumulative hazard and its greatest convex minorant is of the order $o_{p}(n^{-1/2})$,
again with the motivation of establishing asymptotic optimality of the constrained estimator.
A similar result is proved in Kochar, Mukerjee and Samaniego~\cite{kochar2000estimation} for a monotone mean residual life function.
In the regression setting with a fixed design, Durot and Toquet~\cite{durottocquet2003} consider the supremum distance between the partial sum process
and its least concave majorant and prove that,
if the regression function is decreasing with a continuous derivative that  is bounded and bounded away from zero,
then this supremum distance is of the order $O_{p}(n^{-2/3}(\log n)^{2/3})$.
They also provide a lower bound, showing that $n^{-2/3}(\log n)^{2/3}$ is the exact order of the supremum distance.
A generalization to the case of  a random design was developed by Pal and Woodroofe~\cite{pal2006distance}.
Similar results were proved for other shape-constrained
estimators, see Balabdaoui and Wellner~\cite{balabdaoui-wellner2007} for convex densities
and D\"umbgen and Rufibach~\cite{duembgen2009maximum} for log-concave densities.
Wang and Woodroofe~\cite{wang2007kiefer} obtained a similar result for Wicksell's problem.
Their result compares to the supremum distance between the primitive of~$F_n$ and its least concave majorant,
which leads to a faster rate $n^{-1}\log n$.

Although the first motivation for Kiefer-Wolfowitz type of results has been asymptotic optimality of shape constrained estimators,
other important statistical applications are conceivable.
For instance, the Kiefer-Wolfowitz result was a key argument in Sen, Banerjee and Woodroofe~\cite{sen2010inconsistency} to prove that,
although bootstrapping from the empirical distribution function $F_{n}$ or from its least concave majorant~$\widehat F_{n}$
does not work for the Grenander estimator of a decreasing density function at a fixed point,
the $m$ out of $n$ bootstrap, with $m\ll n$, from~$\hat F_{n}$ does work.
Likewise, Durot, Groeneboom and Lopuha\"a~\cite{durot-groeneboom-lopuhaa2013} use a Kiefer-Wolfowitz type of result to prove that a smoothed bootstrap from
a Grenander-type estimator works for $k$-sample tests in a general statistical setting,
which covers the monotone regression model and monotone density model among others.
Mammen~\cite{mammen1991estimating} suggests to use such a result to make an asymptotic comparison of two different estimators for a monotone regression function:
one of them is obtained by smoothing a Grenander type estimator and the other one is obtained by ``monotonizing'' a kernel estimator.
See also Wang and Woodroofe~\cite{wang2007kiefer} for a similar application of their Kiefer-Wolfowitz comparison theorem.

The aim of this paper is to establish a Kiefer-Wolfowitz type of result in
a very general setting that covers the setting considered in~\cite{durot-groeneboom-lopuhaa2013}.
We recover the aforementioned Kiefer-Wolfowitz type of results for $\widehat{F}_n-F_n$ as special cases of our general result.
Furthermore, in a similar general setting we consider the supremum distance between the primitive of $F_n$ and its least concave majorant,
and obtain the same faster rate as found in~\cite{wang2007kiefer}.
As an application of our results,
we consider the problem of estimating a smooth monotone function and provide an asymptotic comparison between an ordinary kernel estimator and a smooth monotone estimator.

The paper is organized as follows.
In Section~\ref{sec:main}, we define our general setting and state our Kiefer-Wolfowitz type inequality.
Section~\ref{sec:examples} is devoted to specific settings to which our main theorem applies.
Applications to estimating smooth monotone functions are described in Section~\ref{sec:application}.
Proofs are deferred to Section~\ref{sec:proofs}.

\section{A Kiefer-Wolfowitz type of inequality in a general setting}
\label{sec:main}
First, we define our general setting as well as the notation that will be used throughout the paper.
Then we state our main result.
The result will be illustrated for several classical settings, such as monotone density or monotone regression, in Section \ref{sec:examples}.

\subsection{The setting}\label{sec: setting}
Suppose that based on $n$ observations, we have at hand a cadlag step estimator~$F_{n}$ for a concave function $F:[a,b]\to\R$,
where $a$ and $b$ are know reals.
In the sequel, we assume that~$F$ is continuously differentiable with $F(a)=0$
and we denote by $f$ the first derivative, which means that
\begin{equation}
\label{eq:def F}
F(t)=\int_{a}^tf(x)\, \mathrm{d}x,
\end{equation}
for $t\in[a,b]$.
A typical example is the case where we have independent observations with a common density~${f}$ on~$[a,b]$,
and where the estimator for $F$ is the empirical distribution function $F_n$ of the observations.
Further details are given in Section~\ref{sec:examples}, where some more examples are investigated.

We will impose the following assumptions on $f$:
\begin{enumerate}
\item[\mono]
The function $f:[a,b]\mapsto \mathbb{R}$ is decreasing and continuously differentiable,
such that $0<\inf_{t\in[a,b]}|f'(t)|\leq \sup_{t\in[a,b]}|f'(t)|<\infty$.
\end{enumerate}
Furthermore, we assume that the cadlag estimator $F_n$ can be approximated in the sense that
\begin{equation}
\label{eq:embedGen}
\sup_{t\in[a,b]}|F_{n}(t)-F(t)-n^{-1/2}B_{n}\circ L(t)|
=
O_{p}\left(\gamma_n\right),
\end{equation}
where $\gamma_n\to 0$,  $L:[a,b]\to\R$ is non-decreasing, and $B_{n}$ is a process on  $[L(a),L(b)]$ that satisfies the following two conditions for a given $\tau\in[0,4)$:
\begin{enumerate}
\item[\Bn]
There are positive $K_1,K_{2}$ such that for all $x\in [L(a),L(b)]$, $u\in(0,1]$, and~$v>0$,
$$
\prob\left(
\sup_{|x-y|\leq u}|B_{n}(x)-B_{n}(y)|>v
\right)
\leq
K_{1} \exp(-K_{2}v^2u^{-\tau}).
$$
\item[\Bnn]
There are positive $K_1,K_{2}$ such that for all $x\in [L(a),L(b)]$, $u\in(0,1]$, and~$v>0$,
$$
\prob\left(\sup_{z\geq u}\left\{B_{n}(x-z)-B_{n}(x)-v z^2\right\}>0\right)\leq K_{1}\exp\left(-K_{2}v^2u^{4-\tau}\right).
$$
\end{enumerate}
Finally, we will impose the following smoothness condition on $L$.
\begin{enumerate}
\item[\regul]
The function $L:[a,b]\mapsto \mathbb{R}$ is increasing and continuously differentiable,
such that $0<\inf_{t\in[a,b]}L'(t)\leq \sup_{t\in[a,b]}L'(t)<\infty$.
\end{enumerate}
A typical example is estimation of a monotone density $f$,
in which case the empirical distribution function $F_n$ can be approximated by a sequence of Brownian bridges $B_n$,
$L$ is equal to the cumulative distribution function  $F$ corresponding to $f$, and $\gamma_n=(\log n)/n$ in~\eqref{eq:embedGen},
due to the Hungarian embedding (see~\cite{komlosmajortusnady1975}).
Other examples are the monotone regression model and the random censorship model with a monotone hazard,
in which case $B_n$ is a Brownian motion and $\gamma_n$ relies on the integrability of the errors in the regression case, and  $\gamma_n=(\log n)/n$ in the random censorship model,
see Sections~\ref{sec:reg} and~\ref{subsec:monotone hazard} for more details.

\subsection{Main results}
\label{subsec:main resuls}
Hereafter, $\hat F_n$ denotes the least concave majorant of $F_{n}$ on $[a,b]$.
We are interested in the supremum distance between $F_{n}$ and $\hat F_{n}$.
Our main result is a Kiefer-Wolfowitz type of inequality for the supremum distance in our general setting.
We will obtain such an inequality by decomposing $\hat F_n-F_n$ into two parts,
the difference between the approximating process $F_{n}^B=F+n^{-1/2}B_{n}\circ L$
and its least concave majorant, and remainder terms that can be bounded by means of~\eqref{eq:embedGen}.
We then first establish a Kiefer-Wolfowitz type of result for $F_{n}^B$ by making use of assumptions~\Bn\ and~\Bnn.
The reason is that in typical examples,
the bound provided by the approximation in~\eqref{eq:embedGen} is of smaller order than
the bound on the difference between $F_n^B$ and its least concave majorant,
and hence, the latter difference determines the rate in the
Kiefer-Wolfowitz result.
Moreover, it has the advantage that it allows one to avoid the specific structure
of the particular statistical model at hand, and it only requires assumptions~\Bn\ and~{\Bnn}
on the approximating process.
Note however, that if the specific structure does provide suitable exponential bounds on
tail probabilities for $n^{1/2}(F_n-F)$,
one can just take the identity for $L$ and $B_n=n^{1/2}(F_n-F)$ in~\eqref{eq:embedGen}.
As it may be of interest in its own right, we first state a Kiefer-Wolfowitz type of result for $F_{n}^B$.
\begin{theorem}
\label{theorem:KWGenB}
Let $F_{n}^B=F+n^{-1/2}B_{n}\circ L$, where $F$ is defined by \eqref{eq:def F} for some~$f$ satisfying~{\mono},
$L$ satisfying {\regul}, and $B_{n}$ satisfying {\Bn} and {\Bnn} for some $\tau\in[0,4)$.
Let $\hat F_{n}^B$ be the least concave majorant of~$F_{n}^B$ on~$[a,b]$.
We then have
$$
\sup_{x\in[a,b]}|\hat F_n^B(x)-F_{n}^{B}(x)|=O_{p}\left(\frac{\log n}n\right)^{2/(4-\tau)}.
$$
\end{theorem}
The main ingredient to prove Theorem~\ref{theorem:KWGenB} is a localization result stated in
Lemma~\ref{lem:espGunifGen} below.
It shows that
although the least concave majorant $\hat F_{n}^B$ depends on the whole process~$F_{n}^B$,
its value at a fixed point $x$ mainly depends on  $F_{n}^B$ in a small neighborhood  of~$x$.
Precisely, with probability tending to one,
$\hat F_{n}^B(x)$  coincides with  the least concave majorant of the restriction of $F_{n}$ to a shrinking interval with center~$x$.
This result generalizes Lemma~5.1 in \cite{durottocquet2003}, where only the case of a Brownian motion $B_{n}$ with the specific variance function $L(t)=t$ was considered.
\begin{lemma}
\label{lem:espGunifGen}
Assume the conditions of Theorem \ref{theorem:KWGenB}. Let
\begin{equation}\label{eq:defcnGen}
c_{n}=\left(\frac{c_{0}\log n}n\right)^{1/(4-\tau)}
\end{equation}
for some $c_{0}>0$.
For $x\in[a,b]$, let  $\hat F_{n,c_{n}}^{(B,x)}$ be the least concave majorant of the process
$\left\{F_n^B(\eta), \eta\in[x-2c_{n},x+2c_{n}]\cap[a,b]\right\}$.
Then, there exist positive numbers $K_{1},K_{2},C_{0}$ independent of $n$, such that for $c_{0}\geq C_{0}$ we have
$$
\prob\left(\sup_{x\in[a,b]}
\left|\hat F_n^B(x)-\hat F_{n,c_n}^{(B,x)}(x)\right|\neq 0\right)
\leq K_{1} n^{-c_{0}K_{2}}.
$$
\end{lemma}
Theorem~\ref{theorem:KWGenB}, together with~\eqref{eq:embedGen}, yields the following general Kiefer-Wolfowitz type of result
for the cadlag estimator $F_n$ and its least concave majorant $\hat F_n$.
\begin{theorem}
\label{theorem:KWGenE}
Assume~\eqref{eq:embedGen}, where $F$ is defined by~\eqref{eq:def F} for some $f$
satisfying~{\mono}, $L$~satisfying {\regul}, and $B_{n}$ satisfying {\Bn} and {\Bnn} for some $\tau\in[0,4)$.
We then have
$$
\sup_{x\in[a,b]}|\hat F_n(x)-F_{n}(x)|
=
O_p\left(\gamma_n\right)+O_{p}\left(\frac{\log n}n\right)^{2/(4-\tau)}.
$$
\end{theorem}
Obviously, when the approximation in~\eqref{eq:embedGen} is sufficiently strong, that is, $\gamma_n=O(n^{-1}\log n)^{2/(4-\tau)}$ then
\begin{equation}
\label{eq:KWresult}
\sup_{x\in[a,b]}|\hat F_n(x)-F_{n}(x)|
=
O_{p}\left(\frac{\log n}n\right)^{2/(4-\tau)}.
\end{equation}
For models where~{\Bn} and~{\Bnn} hold with $\tau=1$ and $\gamma_n=O(n^{-1}\log n)^{2/3}$ in~\eqref{eq:embedGen},
we recover the traditional Kiefer-Wolfowitz inequality (\cite{kiefer-wolfowitz1976}).
See Section~\ref{sec:examples} for examples.
For models where~{\Bn} and~{\Bnn} hold with $\tau=2$ and $\gamma_n=O(n^{-1}\log n)$, we recover the faster rate found by Wang and Woodroofe (see Theorem 2.1 in~\cite{wang2007kiefer}).
The reason for finding different values for $\tau$, is that the case $\tau=1$
corresponds to a Kiefer-Wolfowitz inequality derived for a naive estimator $F_n$ for $F$ in~\eqref{eq:def F},
whereas the result in~\cite{wang2007kiefer} compares to an inequality for the integral of~$F_n$.
See Section~\ref{sec:positivity} for more details.

Under slightly more restrictive assumptions, the
results in Theorems~\ref{theorem:KWGenB} and~\ref{theorem:KWGenE}
can be made more precise by considering moments
of the supremum distance rather than the stochastic order.
As before, we first obtain a result for moments corresponding to the process $F_n^B$.
\begin{theorem}
\label{theorem:KWSBGenEsp}
Assume the conditions of Theorem \ref{theorem:KWGenB}.
Moreover, assume that there are positive numbers $K_{1},K_{2}$,
such that for all  $v>0$ we have
\begin{equation}
\label{eq:BnSup}
\prob\left(\sup_{x\in [L(a),L(b)]}
|B_{n}(x)|>v\right)\leq K_{1} \exp(-K_{2}v^2).
\end{equation}
With  $r\geq 1$ arbitrary, we then have
\[
\E\Bigg[\sup_{x\in[a,b]}
\left|\hat F_n^B(x)-F_{n}^{B}(x)\right|^r
\Bigg]
=O\left(\frac{\log n}{n}\right)^{2r/(4-\tau)}.
\]
\end{theorem}
A similar result for the process $F_n$ is obtained from the previous theorem and
the following condition
\begin{equation}
\label{eq:moments}
\E
\Bigg[
\sup_{x\in[a,b]}|F_{n}(x)-F(t)-n^{-1/2}B_{n}\circ L(t)|^r
\Bigg]=O\left(\frac{\log n}{n}\right)^{2r/(4-\tau)},
\end{equation}
where $F$ is defined by~\eqref{eq:def F} for some $f$ satisfying~{\mono},
$L$ satisfying~{\regul}, and~$B_{n}$  satisfying~{\Bn} and {\Bnn}.
Note that the slightly more restrictive moment assumption \eqref{eq:moments} replaces condition~\eqref{eq:embedGen},
that was used before in Theorem~\ref{theorem:KWGenE}.
\begin{theorem}
\label{theorem:KWEGenEsp}
Assume the conditions of Theorem \ref{theorem:KWGenB}.
Moreover, assume that~\eqref{eq:moments} holds, for some $r\geq 1$ and $\tau\in[0,4)$,
and assume that there are positive numbers $K_{1},K_{2}$ such that~\eqref{eq:BnSup} holds for all  $v>0$.
We then have
\[
\E\Bigg[\sup_{x\in[a,b]}
\left|\hat F_n(x)-F_{n}(x)\right|^r
\Bigg]
=O\left(\frac{\log n}{n}\right)^{2r/(4-\tau)}.
\]
\end{theorem}

\subsection{Local version}
\label{sec:local}
We also investigate a local version of the Kiefer-Wolfowitz result.
This means that instead of considering the supremum over the whole interval $[a,b]$ as in Theorem~\ref{theorem:KWGenE},
we consider the supremum over a shrinking neighborhood around a fixed point $x_{0}\in[a,b]$.
For the local supremum, we obtain a bound of smaller order than for the global supremum.
This compares to Theorem~2.2 in \cite{wang2007kiefer}, where a Kiefer-Wolfovitz type of result was established for Wicksell's problem.
However, only one specific rate of shrinking was considered in~\cite{wang2007kiefer},
whereas we allow a range of possible rates.
Moreover, we give a more precise bound than in~\cite{wang2007kiefer}.
\begin{theorem}
\label{theorem:localrate}
Fix $x_{0}\in[a,b]$.
Assume~\eqref{eq:embedGen}, where $F$ is defined by~\eqref{eq:def F} for some~$f$
satisfying~{\mono}, $L$~satisfying {\regul}, and $B_{n}$ satisfying {\Bn} and {\Bnn} for some $\tau\in[0,4)$.
For any sequence $\eps_{n}\geq (n^{-1}\log n)^{1/(4-\tau)}$ we then have
\[
\sup_{|x-x_{0}|\leq\eps_{n}}|\hat F_{n}(x)-F_{n}(x)|
=
O_p(\gamma_n)
+
O_{p}\left((\eps_{n}^{\tau/2}n^{-1/2})\\\wedge \left(\frac{\log n}n\right)^{2/(4-\tau)}\right).
\]
\end{theorem}
As in the case of Theorem~\ref{theorem:KWGenE}, when the embedding in~\eqref{eq:embedGen} is sufficiently strong, i.e.,
$\gamma_n=O(\eps_{n}^{\tau/2}n^{-1/2})$ and $\gamma_n=O(n^{-1}\log n)^{2/(4-\tau)}$, we obtain
\begin{equation}
\label{eq:theolocal}
\sup_{|x-x_{0}|\leq\eps_{n}}|\hat F_{n}(x)-F_{n}(x)|
=
O_{p}\left((\eps_{n}^{\tau/2}n^{-1/2})\\\wedge \left(\frac{\log n}n\right)^{2/(4-\tau)}\right).
\end{equation}
Clearly, the local rate  in~\eqref{eq:theolocal} is at most $(n^{-1}\log n)^{2/(4-\tau)}$
and for any allowable sequence $\eps_{n}\geq  (n^{-1}\log n)^{1/(4-\tau)}$, it is at least
$(n^{-1}\log n)^{2/(4-\tau)}(\log n)^{-1/2}$.
Thus, the local rate may vary depending on the rate $\eps_n$ at which the neighborhood around $x_0$ shrinks,
and it is of smaller order than the global rate obtained in Theorem~\ref{theorem:KWGenE}
in all cases where $\eps_{n}=o\left((n^{-1}\log n)^{1/(4-\tau)}(\log n)^{1/\tau}\right)$.
\begin{remark}
\label{rem:wicksell local}
Note that for $\tau=2$, the boundary case $\eps_{n}=(n^{-1}\log n)^{1/2}$ coincides with the shrinking rate
in Theorem~2.2 in~\cite{wang2007kiefer}.
This leads to local rate $O_p(\eps_n^2(\log n)^{-1/2})$ in~\eqref{eq:theolocal},
which is conform the the rate $o_p(\eps_n^2)$, as stated in Theorem~2.2
in~\cite{wang2007kiefer}.
\end{remark}
We end this section by considering the rate of convergence at a fixed point.
As stated in the following theorem, the resulting rate is $n^{-2/(4-\tau)}$ with no logarithmic term.
\begin{theorem}\label{lem:pointwiserate}
Fix $x_{0}\in[a,b]$ and suppose that $\gamma_n=O(n^{-2/(4-\tau)})$  in~\eqref{eq:embedGen}.
Then under the assumptions of Theorem \ref{theorem:localrate} we have
\begin{equation}\label{eq:pointwiserate}
\hat F_{n}(x_{0})-F_{n}(x_{0})=O_{p}\left(n^{-2/(4-\tau)}\right).
\end{equation}
\end{theorem}
For models where~{\Bn} and~{\Bnn} hold with $\tau=1$, such as the monotone density model, the rate $n^{-2/3}$ in
Theorem~\ref{lem:pointwiserate} matches with the result in~\cite{wang1994}.
See also~\cite{durottocquet2003} and~\cite{kulikov-lopuhaa2006SPL}.

\section{Examples of specific settings}
\label{sec:examples}
The section is devoted to specific settings to which Theorem~\ref{theorem:KWGenE} applies.
We first discuss statistical models for which a Kiefer-Wolfowitz result is obtained
for an estimator $F_n$ for the integral $F$ of a decreasing curve $f$, and in which the approximation
in~\eqref{eq:embedGen} is by means of Brownian motion of Brownian bridge.
In these cases the Kiefer-Wolfowitz result coincides with the traditional one in~\cite{kiefer-wolfowitz1976}.
Next, we consider the situation for which a Kiefer-Wolfowitz result is obtained for the primitive of $F_n$.
This matches the setup for the Wicksell problem considered by~\cite{wang2007kiefer} and we obtain the
same (faster) rate as found in~\cite{wang2007kiefer}. Finally, we discuss a few setups that are not covered by our general setting.

\subsection{Decreasing functions}
\label{sec:monoGen}
It turns out that in various usual settings (where the decreasing function~$f$ could be for instance a density, or a regression function), the embedding \eqref{eq:embedGen} holds with $B_{n}$ being either Brownian Bridge or Brownian motion.  For such a process $B_{n}$, it can be proved that {\Bn} and {\Bnn} hold with $\tau=1$, leading to the usual rate $(n^{-1}\log n)^{2/3}$ in the Kieffer-Wolfowitz inequality. This is made precise in the following corollary.
Then, we discuss a number of specific settings that are covered by Corollary~\ref{cor:GenMono}.
\begin{corollary}
\label{cor:GenMono}
Assume \eqref{eq:embedGen} with $\gamma_n=O(n^{-1}\log n)^{2/3}$,
where $F$ is defined by~\eqref{eq:def F} for some $f$ that satisfies {\mono}, $L$ satisfies {\regul}, and $B_{n}$  is either Brownian motion or Brownian Bridge. We then have
\begin{equation}\label{eq:rate2/3}
\sup_{x\in[a,b]}
\left|\hat F_n(x)-F_{n}(x)\right|=O_{p}\left(\frac{\log n}n\right)^{2/3}.
\end{equation}
If, moreover, $0\in[L(a),L(b)]$ and \eqref{eq:moments} holds with $\tau=1$ and some $r\geq 1$, then we also have
\begin{equation}\label{eq:rate2/3moments}
\E\left[\sup_{x\in[a,b]}
\left|\hat F_n(x)-F_{n}(x)\right|^r\right]=O\left(\frac{\log n}n\right)^{2r/3}.
\end{equation}
\end{corollary}

\subsubsection{Monotone regression function.}\label{sec:reg}
We have observations $Y_{i}$, for $i=1,2,\ldots,n$, satisfying
$Y_{i}={f}(t_{i})+\epsilon_{i}$,
where $\E(\epsilon_{i})=0$ and
\begin{equation}\label{eq:design}
\max_{i}|t_{i}-(a+(b-a)i/n)|=O(n^{-2/3}),
\end{equation}
 which means that the design points are close to uniformly spread on $[a,b]$.
We assume that the $\epsilon_{i}$'s are independent having the same distribution.
In this case,  the estimator for~$F$ in~\eqref{eq:def F} is
\begin{equation}
\label{eq:def Fn regression}
F_{n}(t)
=
\frac{1}{n}\sum_{i=1}^{n} Y_{i}\indicator\{t_{i}\leq t\},\ t\in[a,b].
\end{equation}
As a special case of  Corollary \ref{cor:GenMono} we obtain the following result.
\begin{corollary}
\label{cor:regression}
If $\E|\epsilon_{i}|^3<\infty$, $\E|\epsilon_{i}|^2>0$
and  {\mono} holds, then we have \eqref{eq:rate2/3}  with $F_n$ taken from~\eqref{eq:def Fn regression}
and $\hat F_{n}$ the least concave majorant of $F_{n}$.
\end{corollary}
\subsubsection{Monotone density.}
We have independent observations ${X}_{i}$, for $i=1,2,\ldots,n$,
with common density~${f}:[a,b]\to\R$, where $a$ and $b$ are known real numbers.
The estimator for the distribution function $F$ in this case is the empirical distribution function
\begin{equation}
\label{eq:def Fn density}
F_{n}(t)
=
\frac{1}{n}\sum_{i=1}^{n} \indicator\{X_{i}\leq t\},\ t\in[a,b].
\end{equation}
\begin{corollary}
\label{cor:density}
If  \mono~holds and  $\inf_{t\in[0,1]}f(t)>0$, then we have \eqref{eq:rate2/3}  with $F_n$ taken from~\eqref{eq:def Fn density}
and $\hat F_{n}$ the least concave majorant of $F_{n}$.\end{corollary}
\subsubsection{Random censorship with monotone hazard.}
\label{subsec:monotone hazard}
We have right-censored observations
$(X_{i},\Delta_{i})$, for $i=1,2,\ldots,n$,
where $X_{i}=\min(T_{i},Y_{i})$ and $\Delta_{i}=\indicator\{T_{i}\leq Y_{i}\}$.
The failure times~$T_{i}$ are assumed to be nonnegative independent with distribution function $G$
and are independent of the i.i.d.~censoring times~$Y_{i}$ that have distribution function $H$.
Define $F=-\log(1-G)$ the cumulative hazard on $[0,b]$.
Note that in this setting, we only consider the case $a=0$, since this is more natural.
The estimator for the cumulative hazard $F$ is defined via the Nelson-Aalen estimator~$N_{n}$ as follows:
let $t_{1}<\cdots<t_{m}$ denote the ordered distinct uncensored failure times in the  sample
and $n_{k}$ the number of $i\in\{1,2,\dots,n\}$ with $X_{i}\geq t_{k}$,\begin{equation}
\label{eq:def Fn censoring}
F_{n}(t_{i})
=
\sum_{k\leq i}\frac{1}{n_{k}},
\end{equation}
and $F_{n}(t)=0$ for all $t<t_{1}$ and $F_{n}(t)=N_{n}(t_{m})$ for all $t\geq t_{m}$.
\begin{corollary}
\label{cor:censoring}
Suppose \mono,
$\inf_{t\in[0,b]}f(t)>0$,
$G(b)<1$, and $\lim_{t\uparrow b}H(t)<1$.
Then we have \eqref{eq:rate2/3}  with~$F_n$ taken from~\eqref{eq:def Fn censoring}
and
$\hat F_{n}$ the least concave majorant of the restriction of $F_{n}$ to $[0,b]$.
\end{corollary}

\subsection{Decreasing primitive of nonnegative functions}
\label{sec:positivity}
Wang and Woodroofe~\cite{wang2007kiefer} obtain a Kiefer-Wolfowitz result for the Wicksell problem
and find $n^{-1}\log n$ as rate of convergence, which is faster than $(n^{-1}\log n)^{2/3}$ from Corollary~\ref{cor:GenMono}.
The reason is that in their setup the Kiefer-Wolfowitz result is obtained for $U_n^\#$, defined as the primitive
of $\Psi_n^\#$, which serves as an estimator for $\Psi(y)=\int_{y}^\infty \pi^2\varphi(x)\,\mathrm{d}x$, for some nonnegative $\varphi$.
We investigate a similar setup, where we establish a Kiefer-Wolfowitz result for the primitive
of the cadlag estimator $F_n$ for $F$ with $f$ being positive.
Precisely, in the sequel we assume
\begin{enumerate}
\item[\positive]
The function $f:[a,b]\mapsto \mathbb{R}$ is continuous,
such that $0<\inf_{t\in[a,b]}f(t)\leq \sup_{t\in[a,b]}f(t)<\infty$.
\end{enumerate}

Let $F$ be defined in~\eqref{eq:def F} with corresponding cadlag estimator $F_n$.
Define
\begin{equation}
\label{eq:def S}
\begin{split}
S(y)
&=
\int_y^b f(x)\,\mathrm{d}x=F(b)-F(y);\\
S_n(y)
&=
F_n(b)-F_n(y).
\end{split}
\end{equation}
If $f$ satisfies~\positive, then $S$ is decreasing.
We will investigate a Kiefer-Wolfowitz result for
\begin{equation}
\label{eq:def Hn}
H_{n}(t)=\int_{a}^{t} S_{n}(x)\,\mathrm{d}x=\int_{a}^{t} (F_{n}(b)-F_n(x))\,\mathrm{d}x,
\end{equation}
which serves as an estimator for
\begin{equation}
\label{eq:def H}
H(t)=\int_{a}^{t} S(x)\,\mathrm{d}x=\int_{a}^{t} (F(b)-F(x))\,\mathrm{d}x.
\end{equation}
In the case where~\eqref{eq:embedGen} holds with $\gamma_n=O(n^{-1}\log n)$ and $B_n$ is either Brownian motion
or Brownian bridge, we find the same rate as in~\cite{wang2007kiefer}.
Let $\hat H_n$ denote the least concave majorant of $H_n$ on $[a,b]$.

\begin{corollary}
\label{cor:GenPos}
Assume \eqref{eq:embedGen} with $\gamma_n=O(n^{-1}\log n)$,
where $F$ is defined by~\eqref{eq:def F} for some $f$ that satisfies~\positive,
$L$ is bounded, and $B_{n}$  is either Brownian motion or Brownian Bridge.
Then, we have
$$
\sup_{x\in[a,b]}|\hat H_n(x)-H_{n}(x)|=O_{p}\left(\frac{\log n}n\right),
$$
\end{corollary}
Note that we recover the rate obtained in~\cite{wang2007kiefer} for Wicksell's problem.
Our result applies for instance to the regression setting of Subsection~\ref{sec:reg},
where the $\eps_{i}$'s are i.i.d.~Gaussian
and, where instead of~{\mono}, $f$ is assumed to satisfy~{\positive}.

\begin{remark}
\label{rem:wicksell}
The general setup of Corollary~\ref{cor:GenPos} does not cover the Wickell problem
considered in~\cite{wang2007kiefer}.
The reason is that the approximating process for their process $U_n^\#$
is more complicated than the one for $H_n$, leading to extra logarithmic terms in~{\Bn} and~{\Bnn}, and to difficulties in obtaining bounds for a supremum in these assumptions.
Nevertheless, by using the specific structure of the Wicksell problem, the authors do obtain the same rate $n^{-1}\log n$,
see Theorem~2.2 in~\cite{wang2007kiefer}.
\end{remark}

\subsection{Discussion}
Although, the above Kiefer-Wolfowitz results have been obtained for two different general settings, 
there are still complex statistical models that are not covered by either setup.
One such example is interval censoring, where one would expect the same result as in Theorem~\ref{theorem:KWGenE}.
The main difference with our general setting is that the estimator for $f$ is the slope of the least concave majorant
of a cumulative sum diagram (CSD), which has $y$-coordinates determined by a cadlag function, e.g., $F_n(t)$ 
in the case of monotone density or monotone regression, and an $x$-coordinate determined by the identity.
In this case, the points are of the form $(t,F_n(t))$ for specific values of $t$,
e.g., $t=X_{(1)},\ldots,X_{(n)}$ in the case of monotone density.
This is essentially different from, for instance, interval censoring Case 1.
In this model, the observations are $(T_i,\Delta_i)$, where $\Delta_i=\{X_i\leq T_i\}$ indicates
whether the event time $X_i$ is before or after the observed censoring time $T_{i}$.
The parameter of interest is the cdf $F$ of the event times
and the coordinates of the CSD are of the form
\[
\left(
G_n(s),\int\{x\leq t\}\{t\leq s\}\,\mathrm{d}\mathbb{P}_n(x,t)
\right),
\quad
s=T_{(1)},\ldots,T_{(n)},
\]
where $G_n$ is the empirical cdf of the $T_i$ and $\mathbb{P}_n$ is the empirical measure of the tuples~$(X_i,T_i)$.
First of all, one would need to find an approximating process satisfying (A2)-(A3), for the process
\[
Y_n(s)=\int\{x\leq t\}\{t\leq s\}\,\mathrm{d}\mathbb{P}_n(x,t).
\]
More importantly, our proof of the key Lemma~2.1 relies heavily on the fact that the points of the CSD
are of the form $(t,Y_n(t))$, for some process $Y_n$, and it is not clear how this can be extended to a CSD with points
of the form $(G_n(t),Y_n(t))$.
Interval censoring case 2 is even more difficult, because the CSD is self-induced,
i.e., the points of the diagram depend on the actual solution itself.

\section{Estimating a smooth monotone function}
\label{sec:application}
In many applications, the parameter of interest $f:[a,b]\to\R$,
e.g., a density function, a regression mean, or a failure rate,
is known to be non-increasing (the non-decreasing case can be treated likewise)
 so it is natural to incorporate this shape constraint into the estimation procedure.
Consider the setting of Section~\ref{sec:monoGen}.
A popular estimator for $f$ under the constraint that $f$ is non-increasing is the Grenander-type estimator $\hat f_{n}$,
defined on $(a,b]$ as the left-hand slope of the least concave majorant $\hat F_{n}$ of $F_{n}$, with
$$
\hat f_{n}(a)=\lim_{s\downarrow a}\hat f_{n}(s).
$$
This estimator is a step function and as a consequence it is not smooth.
Moreover, the rate of convergence of $\hat f_{n}$ is $n^{1/3}$, if $f$ has a first derivative $f'$ that is bounded away from zero,
whereas competitive smooth estimators may have faster rates in cases where $f$ is smooth.
On the other hand, such estimators typically do not satisfy the monotonicity constraint.

In this section, we are interested in an estimator that is both non-increasing and smooth, and that achieves the optimal rate of convergence under certain smoothness conditions.
The estimator is obtained by smoothing the Grenander-type estimator $\hat f_{n}$,
and resembles the estimators $m_{IS}$ in \cite{mammen1991estimating} and $m_{n}$ in~\cite{mukerjee1988monotone},
see also~$\tilde{\Psi}_{n,s}$ in~\cite{wang2007kiefer}.
In this way, one first applies an isotonization procedure followed by smoothing.
A natural alternative would be to interchange the two steps,
that is, first smooth and then isotonize, but this typically results in a non smooth estimator.
It may happen that the two proposals are asymptotically equivalent in first order;
see \cite{mammen1991estimating} for a precise statement in the smooth regression setting.
See also~\cite{vandervaart-vanderlaan2003} for a comparison of the second proposal with an ordinary kernel estimator and with the Grenander estimator when estimating a monotone density with a single derivative.

Consider an ordinary kernel-type estimator~$\widetilde f_n$, corrected at the boundaries in such a way that it converges to $f$,
with a fast rate over the whole interval~$[a,b]$ (whereas the non-corrected kernel estimator may show difficulties at the boundaries):
for every $t\in[a+h_{n},b-h_{n}]$,
\begin{equation}
\label{eq:tildefmain}
\widetilde f_n(t)
=
\frac{1}{h_n}\int_\R K\left(\frac{t-x}{h_n}\right)\,\md F_{n}(x),
\end{equation}
where $h_n>0$ and the kernel function $K:\R\to[0,\infty)$ satisfies~$\int K(t)\,\md t=1$. We are interested in $\widehat f_{ns}$, the estimator defined in the same manner as $\widetilde f_{n}$, with~$F_{n}$ replaced by the least concave majorant~$\widehat{F}_{n}$.
At the boundaries $[a,a+h_{n})$ and $(b-h_{n},b]$, we consider the local linear bias correction defined as follows:
with~$f_{n}$ denoting either $\widetilde f_{n}$ or $\hat f_{ns}$,
\begin{equation}
\label{eq:tildefKosorok}
f_n(t)
=
\begin{cases}
f_n(a+h_n)+ f_n'(a+h_n)(t-a-h_n), &t\in[a,a+h_n];\\
\\[-10pt]
f_n(b-h_n)+f_n'(b-h_n)(t-b+h_n), &t\in[b-h_n,b],
\end{cases}
\end{equation}
see, e.g.~\cite{wand-jones1995}.
Thus, $\hat f_{ns}$ is a smoothed version of the Grenander-type estimator~$\hat f_{n}$, linearly extended at the boundaries.
According to the following lemma, it is monotone provided that $K\geq 0$ is supported on $[-1,1]$.
A similar result was obtained by \cite{mukerjee1988monotone}, page 743, in the regression setting for  a log-concave kernel $K$.
Moreover, since $\hat f_{n}$ can easily be computed using the PAVA or a similar device, see e.g., \cite{barlow1972statistical},
the monotone smooth estimator $\widehat f_{ns}(t)$ is easy to implement thanks to~\eqref{eq:implement} below.
This was already pointed out in~\cite{durot-groeneboom-lopuhaa2013}, Section~4.2.
\begin{lemma}
\label{lem:properties f_ns}
Let $p_1,\dots,p_m$ be the jump sizes of  $\hat f_{n}$ at the points of jump $\tau_1<\dots<\tau_m\in(a,b]$. If $K\geq 0$ is supported on $[-1,1]$, then $\hat f_{ns}$ is non-increasing on~$[a,b]$ and for all $t\in[a+h_{n},b-h_{n}]$, we have
\begin{equation}
\label{eq:implement}
\hat f_{ns}(t)
=\sum_{j=1}^m p_j\int_{(t-\tau_j)/h_{n}}^{\infty}K(u)\,\md u+\hat f_{n}(b).
\end{equation}
\end{lemma}
As application of Corollary~\ref{cor:GenMono}, we establish that $\widehat{f}_{ns}$ is uniformly
close to~$\widetilde f_{n}$,
and similarly for their derivatives.
This will ensure that the two estimators and their derivatives are asymptotically equivalent in first order.
In~\cite{wang2007kiefer} a similar application of a Kiefer-Wolfowitz result is discussed.
Their result is for $t$ fixed and compares to our result for the derivatives of $\widehat{f}_{ns}$ and $\widetilde{f}_n$.
\begin{lemma}
\label{lem:compar}
Assume the conditions of Corollary~\ref{cor:GenMono}. If $K$ is supported on $[-1,1]$ with integrable first and second derivatives,
then for $l=0,1$,
$$
\sup_{t\in[a,b]}\left|
\widehat{f}_{ns}^{(l)}(t)-\widetilde{f}_n^{(l)}(t)
\right|=O_{p}\left(h_{n}^{-(1+l)}n^{-2/3}(\log n)^{2/3}\right),
$$
where $\widehat{f}_{ns}^{(l)}$ and $\widetilde{f}_n^{(l)}$ denote $l$-th derivatives.
\end{lemma}
Thanks to Lemma~\ref{lem:compar}, we are able to derive the limit behavior of $\widehat f_{ns}$ from that of $\widetilde f_{n}$.
To illustrate this, suppose that $f$ belongs to a H\"older class~${\cal H}(L,\alpha)$,
for some $L>0$ and $\alpha\in (1,2]$, which means that $f$ has a first derivative satisfying
$$
|f'(u)-f'(v)|\leq L|u-v|^{\alpha-1},
$$
for all $u,v\in[a,b]$. It is known that in typical settings (including the specific settings investigated in Subsection~\ref{sec:monoGen}),
the estimator defined by~\eqref{eq:tildefmain} with
\begin{equation}
\label{eq:hn}
h_{n}=
R_{n}n^{-1/(2\alpha+1)},
\end{equation}
where $0<R_n+R_{n}^{-1}=O_{P}(1)$, and a kernel function $K$ with
$\int uK(u)\,\md u=0$,
satisfies
\[
\widetilde f_{n}(x)-f(x)=O_{P}\left(n^{-\alpha/(2\alpha+1)}\right),
\]
for all fixed $x\in(a,b)$ independent of $n$. Moreover,  this rate of convergence is optimal in the minimax sense
in typical settings, e.g., see Theorem~2.3 in \cite{cybakov2003introduction}.
With $h_{n}$ defined as in~\eqref{eq:hn}, Lemma \ref{lem:compar} yields  that
\[
\widehat f_{ns}(x)-f(x)
=
\widetilde f_{n}(x)-f(x)
+
o_{P}\left(n^{-\alpha/(2\alpha+1)}\right).
\]
This means that $\widehat f_{ns}$ is asymptotically equivalent to $\widetilde f_{n}$ in first order.
In particular, $\widehat f_{ns}(x)$ has the same limit distribution and the same minimax rate of convergence as $\widetilde f_{n}(x)$,
provided that $h_{n}$ is chosen according to~\eqref{eq:hn}.
Therefore, one can use any adaptive method for calibrating the bandwidth $h_{n}$ of the ordinary kernel estimator $\widetilde f_{n}$, e.g., see~\cite{lepski1997optimal}, and use the same bandwidth in~$\widehat f_{ns}(x)$,
so that it achieves the minimax rate.
Similar arguments enable us to derive the global limit behavior of $\widehat f_{ns}$ from that of $\widetilde f_{n}$,
e.g., the limit distribution or the rate of convergence of the supremum distance between $\widehat f_{ns}$ and $f$.
See~\cite{durot-groeneboom-lopuhaa2013} for further details.

\section{Proofs}
\label{sec:proofs}

Note that it suffices to prove the results for the case $[a,b]=[0,1]$.
Indeed, suppose that $f(t)$, for $t\in[a,b]$, satisfies conditions~\mono,  \eqref{eq:embedGen} and~\regul\ with
corresponding $F$, $L$ and $F_{n}$ on $[a,b]$.
Then this case can be transformed to the case~$[0,1]$ by
considering $(b-a)f(a+x(b-a))$, for $x\in[0,1]$.
It is straightforward to see that these are functions on $[0,1]$ that
satisfy~\mono,  \eqref{eq:embedGen} and~\regul\ with
corresponding functions $F(a+x(b-a))$, $L(a+x(b-a))$ and $F_{n}(a+x(b-a))$ for $x\in[0,1]$.
Moreover, note that the transformed estimator $\hat F_{n}(a+x(b-a))$ is the least concave majorant of the process
$\{F_{n}(a+u(b-a)), u\in[0,1]\}$ at the point $u=x$.
Thus, we prove the results only for the case $[a,b]=[0,1]$.

In the remainder of the section, we assume that $c_{0}$ in \eqref{eq:defcnGen} satisfies $c_{0}\geq C_{0}$ for a given $C_{0}$. The letters $K_{1},K_{2},\ldots$ denote positive numbers that depend only on $f$, $L$ and $C_{0}$ and that may change from line to line. Moreover, we denote $\eps=\inf_{t}|f'(t)|$.

\subsection{Proofs for Subsection \ref{subsec:main resuls}}
Before establishing the key Lemma~\ref{lem:espGunifGen}, we obtain the following result on the increments of $B_{n}$,
which will be used several times.
\begin{lemma}
\label{lem:increments WGen}
Let $B_{n}$ be a process that satisfies~{\Bn} on an interval $I$.
Then there exist  positive $K_1,K_2$ such that for all $u\in(0, 1/2]$ and $v>0$,\[
\prob
\left(
\sup_{x\in I}\sup_{|x-y|\leq u}
\left|
B_{n}(x)-B_{n}(y)
\right|
>
v
\right)
\leq
K_1u^{-1}
\exp
\left(
-K_2v^2u^{-\tau}
\right).
\]
\end{lemma}
\paragraph{Proof.}
Denote $I=[\alpha,\beta]$, let $k$ be the integer part of $(\beta-\alpha)u^{-1}$
and let $t_j=\alpha+ju$, for $j=0,1,\ldots,k+1$.
We then have
\[
\begin{split}
&\prob
\left(
\sup_{x\in[\alpha,\beta]}\sup_{|x-y|\leq u}
\left|
B_{n}(x)-B_{n}(y)
\right|
>
v
\right)\\
&\quad\quad\leq
\sum_{j=0}^{k}
\prob
\left(
\sup_{x\in[t_j,t_{j+1}]}\sup_{|x-y|\leq u}
\left|
B_{n}(x)-B_{n}(y)
\right|
>
v
\right).
\end{split}
\]
Since $t_{j+1}-t_{j}=u$, for all $j=0,\dots k$ we have
\[
\begin{split}
&\prob
\left(
\sup_{x\in[t_j,t_{j+1}]}\sup_{|x-y|\leq u}
\left|
B_{n}(x)-B_{n}(y)
\right|
>
v
\right)\\
&\quad\quad\leq
\prob
\left(
\sup_{x\in[t_j,t_{j+1}]}
\left|
B_{n}(x)-B_{n}(t_j)
\right|
>
\frac v2
\mbox{ or }\sup_{|t_j-y|\leq 2u}
\left|
B_{n}(t_j)-B_{n}(y)
\right|
>
\frac v2
\right)\\
&\quad\quad\leq
\prob
\left(
\sup_{|t_j-y|\leq 2u}
\left|
B_{n}(t_j)-B_{n}(y)
\right|
>
\frac v2
\right)\\
&\quad\quad\leq
K_{1}\exp(-K_{2}v^2u^{-\tau}/2^{2+\tau}).
\end{split}
\]
We used {\Bn} for the last inequality. We conclude that
\[
\begin{split}
&\prob
\left(
\sup_{x\in[\alpha,\beta]}\sup_{|x-y|\leq u}
\left|
B_{n}(x)-B_{n}(y)
\right|
>
v
\right)\\
&\quad\leq
K_{1}(k+1)\exp\left(-\frac{K_{2}v^2u^{-\tau}}{2^{2+\tau}}\right)
\leq
\frac{K_{1}}u\left(\beta-\alpha+1\right)\exp\left(-\frac{K_{2}v^2u^{-\tau}}{2^{2+\tau}}\right),
\end{split}
\]
since $k\leq (\beta-\alpha)/u$ and $u\in(0,1]$. This proves the lemma by renaming $K_{1}$ and~$K_{2}$.
\tqed

\paragraph{Proof of Lemma \ref{lem:espGunifGen}.}
The proof is inspired by the proof of Lemma~5.1 in~\cite{durottocquet2003}.
Recall that without loss of generality, $[a,b]=[0,1]$.
For all $x\in[0,1]$, let
$$
\tilde x_i=\inf
\left\{u\geq  (x-2c_n)\vee 0\text{, such that }\hat F_n^B(u)= \hat F_{n,c_n}^{(B,x)}(u)\right\},
$$
with the convention that the infimum of an empty set is $(x+2c_n)\wedge1$, and let
$$
\tilde x_s=
\sup
\left\{u\leq  (x+2c_n)\wedge 1\text{, such that }\hat F_n^B(u)= \hat F_{n,c_n}^{(B,x)}(u)\right\},
$$
with the convention that the supremum of an empty set is $(x-2c_n)\vee0$.
If $\hat F_n^B(u)= \hat F_{n,c_n}^{(B,x)}(u)$ for some $u\leq x$,
and $\hat F_n^B(v)= \hat F_{n,c_n}^{(B,x)}(v)$ for some $v\geq x$,
then we must have $\hat F_n^B= \hat F_{n,c_n}^{(B,x)}$ on the whole  interval $[u,v]$.
Therefore, if for some~$x$ we have $\hat F_n^B(x)\neq \hat F_{n,c_n}^{(B,x)}(x)$, then we must have either
$\tilde x_i> x$ or $\tilde x_s< x$.
Moreover, note that if $\tilde x_i> x\geq 0$, then we must have $x-2c_n>0$.
Otherwise, we would have $F_{n,c_n}^{(B,x)}(0)=F_n^B(0)=\hat F_n^B(0)$,
which would mean that $\tilde x_i=0$.
Similarly, if $\tilde x_s< x\leq 1$, then we must have $x+2c_n<1$.
Therefore, it suffices to prove that
there exist positive $K_{1},K_{2}$ such that
\begin{equation}\label{PxiGen}
\prob\left(\tilde x_i>x  ,\text{ for some } x\in [2c_n,1]\right)
\leq K_{1} n^{-K_{2}c_{0}},
\end{equation}
and
\begin{equation}\label{PxsGen}
\prob\left(\tilde x_s<x  ,\text{ for some } x\in[0,1-2c_n]\right)
\leq
K_{1} n^{-K_{2}c_{0}},
\end{equation}
provided that $c_{0}\geq C_{0}$ for some  sufficiently large $C_{0}$.
We will only prove \eqref{PxiGen}, since~\eqref{PxsGen} can be proven with similar arguments.

If $\tilde x_i> x$ for some $x\in[2c_n,1]$,  then
by definition,
$$
\hat F_n^B(u)
\neq
\hat F_{n,c_n}^{(B,x)}(u),
$$
for all $0\leq u\leq x$.
In that case, there exist
$0\leq y\leq x-2c_n$ and $x\leq z\leq (x+2c_{n})\wedge 1$, such
that the line segment joining $(y,F_n^B(y))$ and $(z,F_n^B(z))$ is above $(t,F_n^B(t))$
for all $t\in(y,z)$. In particular, this line segment is above $(x-c_n,F_n^B(x-c_n))$, which
implies that the slope of
the line segment joining $(y,F_n^B(y))$ and $(x-c_n,F_n^B(x-c_n))$ is
smaller than the slope of
the line segment joining $(z,F_n^B(z))$ and $(x-c_n,F_n^B(x-c_n))$. This means that
$$\frac{F_n^B(y)-F_n^B(x-c_n)}{y-x+c_n}<
\frac{F_n^B(z)-F_n^B(x-c_n)}{z-x+c_n}.$$
For any fixed $\alpha\in\R$, this implies that
$$\mbox{ either }\frac{F_n^B(y)-F_n^B(x-c_n)}{y-x+c_n}<\alpha\mbox{ or }\alpha<
\frac{F_n^B(z)-F_n^B(x-c_n)}{z-x+c_n}.$$
In particular with
$\alpha_x=f(x)+c_n|f'(x)|$ we have
\begin{equation}
\label{eq:KWtildeGen}
\prob\left(\tilde x_i>x  \text{ for some } x\in[2c_n,1]\right)\leq \prob_1+\prob_2,
\end{equation}
where
\[\begin{split}
\prob_1
&=
\prob\Big(
\exists x\in[2c_n,1],
\exists y\in[0,x-2c_n]:\\
&\quad\quad\quad\quad F_n^B(y)-F_n^B(x-c_n)>(y-x+c_n)\alpha_x
\Big),
\end{split}
\]
and
\[\begin{split}
\prob_2
&=
\prob\Big(
\exists x\in[2c_n,1],
\exists z\in[x,(x+2c_n)\wedge 1]:\\
&\quad\quad\quad\quad
F_n^B(z)-F_n^B(x-c_n)>(z-x+c_n)\alpha_x\Big).
\end{split}\]
Furthermore, with
$t_x=c_{n}^2 f'(x)/4$, we have $\prob_1\leq \prob_{1,1}+\prob_{1,2}$, where
$$
\prob_{1,1}
=
\prob\left(\exists x\in[2c_n,1]: F_n^B(x)-F_n^B(x-c_n)>c_{n}
\alpha_x+t_x \right)
$$
and
$$
\prob_{1,2}
=
\prob\left(
\exists x\in[2c_n,1],
\exists y\in[0,x-2c_n]:
F_n^B(x)-F_n^B(y)<(x-y)\alpha_x+t_x\right).
$$
We first consider $\prob_{1,1}$.
From~\mono, the derivative $f'$ (which is defined respectively as the right and the left  derivative of $f$ at the boundary points $0$ and $1$) is negative and uniformly continuous on the compact interval $[0,1]$.
Since $c_{n}$ tends to zero, by using Taylor's expansion, we obtain
$$
F(x)-F(x-c_{n})= c_{n}f(x)+\frac{c_{n}^2}{2}\Big (|f'(x)|+o(1)\Big),
$$
where the $o(1)$ term is uniform in $x\in[2c_n,1]$.
Therefore, with $M_{n}^B=F_{n}^B-F$, we obtain
\[
\begin{split}
\prob_{1,1}
&\leq
\prob\left(
\exists x\in[2c_n,1]:
\Big(M_{n}^B(x)-M_{n}^B(x-c_n)\Big)
>
\frac{c_{n}^2}{4}\left(|f'(x)|+o(1)\right)\right)\\
&\leq
\prob\left(
\sup_{x\in[2c_n,1]}\Big(M_{n}^B(x)-M_{n}^B(x-c_n)\Big)>\frac{c_{n}^2}{8}\inf_{t\in[0,1]}|f'(t)|
\right),
\end{split}
\]
provided $n$ is sufficiently large. By definition, $M_{n}^B=n^{-1/2}B_{n}\circ L$. Moreover,
$\left|L(x)-L(x-c_n)\right|\leq c_{n}||L'||_{\infty},
$ where by assumption {\regul},
$$||L'||_{\infty}:=\sup_{t\in[0,1]}L'(t)<\infty.$$
Using Lemma~\ref{lem:increments WGen}, we conclude  that with $\eps=\inf_{t}|f'(t)|>0$ and $I=[L(2c_n),L(1)]$,
\[
\begin{split}
\prob_{1,1}
&\leq
\prob\left(
\sup_{x\in I}\sup_{|x-y|\leq c_n||L'||_{\infty}}
(B_{n}(x)-B_{n}(y))>\frac{c_{n}^2\sqrt n}{8}\eps\right)\\
&\leq
K_{1}||L'||_{\infty}^{-1}c_{n}^{-1}\exp\left(-\frac{K_{2}\eps^2}{64||L'||_{\infty}^\tau}nc_{n}^{4-\tau}\right)\\
&\leq
K_{1}||L'||_{\infty}^{-1}\left(\frac{n}{c_{0}\log n}\right)^{1/(4-\tau)}n^{-K_{2}\eps^2c_{0}/({64||L'||_{\infty}^\tau)}}.
\end{split}
\]
Possibly enlarging $K_{1}$, this proves that for $c_{0}$ sufficiently large and all $n$,
$$
\prob_{1,1}\leq K_{1}n^{-K_{2}\eps^2c_{0}/({65||L'||_{\infty}^\tau)}}.
$$
Renaming $K_{2}$, we conclude that there exist positive numbers $K_{1}$ and $K_{2}$ that depend only on $f$, $L$ and $C_{0}$ such that
\begin{equation}
\label{eq:bound P11Gen}
\prob_{1,1}
\leq K_{1}n^{-K_{2}c_{0}},
\end{equation}
for all $n$, provided that $c_{0}\geq C_{0}$ for some sufficiently large $C_{0}$.

Next, consider $\prob_{1,2}$.
For all
$x\in[2c_n,1]$ and $z\in[1,x/(2c_n)]$,
let $Y_n(x,z)$ be defined by
$$
Y_n(x,z)=
F_n^B(x-2c_nz)-F_n^B(x)
+2c_{n}\alpha_xz+t_x,
$$
so that
\begin{equation}
\label{KWP12Gen}
\prob_{1,2}
=
\prob\left(
\exists x\in[2c_n,1],
\exists z\in[1,x/(2c_n)]:
 Y_n(x,z)>0
\right).
\end{equation}
Let $\eps=\inf_{t\in[0,1]}|f'(t)|$ and let $a$ be a real number with $a\eps> 2\sup_{t\in[0,1]}|f'(t)|$ (which implies that $a\geq 2$).
Moreover, recall that $\alpha_x=f(x)+c_n|f'(x)|$ and $t_x=c_{n}^2 f'(x)/4$.
Now, distinguish between $z\in[1,a]$ and $z\in[a,x/(2c_n)]$.

For all $z\in[a,x/(2c_n)]$, it follows from Taylor's expansion and the definition of $a$  that
\begin{equation}\label{eq:zGeqa}
F(x-2c_nz)-F(x)
+2c_{n}\alpha_{x} z\leq -\eps c_{n}^2z^2.
\end{equation}
Define $A_n=\{(x,z):x\in[2c_n,1],z\in[a,x/(2c_n)]\}$.
From \eqref{eq:zGeqa} we have
\[
\begin{split}
&
\prob\left(\sup_{(x,z)\in A_n}Y_n(x,z)>0\right)
\\
&\quad\leq
\prob\left(
\sup_{(x,z)\in A_n}
\left\{
M_{n}^B(x-2c_nz)-M_{n}^B(x)-\eps c_{n}^2z^2
\right\}>\frac{c_{n}^2\eps}{4}
\right)\\
&\quad\leq
\prob\left(
\sup_{(x,z)\in A_n}
\left\{
B_{n}\circ L(x-2c_nz)-B_{n}\circ L(x)-\eps c_{n}^2\sqrt n z^2
\right\}>\frac{c_{n}^2\eps\sqrt n}{4}
\right).
\end{split}
\]
Define
$A_n'=\{(t,u):t=L(x),u=(L(x)-L(x-2c_{n}z))/(2c_{n}),(x,z)\in A_n\}$.
Then
\[
\begin{split}
&
\prob\left(
\sup_{(x,z)\in A_n}Y_n(x,z)>0\right)
\\
&\leq
\prob\left(
\sup_{(t,u)\in A_n'}
\left\{
B_{n}(t-2c_n u)-B_{n}(t)- \frac{c_{n}^2\eps \sqrt{n}u^2}{||L'||_{\infty}^2}
\right\}>\frac{c_{n}^2\eps\sqrt n}{4}
\right).
\end{split}
\]
Now, denote by $k_n$ the integer part of $c_{n}^{-1}$ and
for all
$j=0,1,\ldots, k_n$, let $t_j=L(2c_n)+j(L(1)-L(2c_n))/k_n$.
If for some $(t,u)\in A_n'$, one has
\[
B_{n}(t-2c_nu)-B_{n}(t)
-
\frac{c_{n}^2\eps \sqrt{n}u^2}{||L'||_{\infty}^2}
>
\frac{c_{n}^2\eps\sqrt n}{4},
\]
then, for $j=1,2,\ldots, k_n$, such that $t\in[t_{j-1},t_{j}]$, one either has
\[
B_{n}(t_{j}-2c_nu)-B_{n}(t_{j})-\frac{c_{n}^2\eps \sqrt{n}u^2}{||L'||_{\infty}^2}
>0,
\]
or
\[
B_{n}(t-2c_nu)-B_{n}(t_{j}-2c_nu)
-
B_{n}(t)+B_{n}(t_{j})>\frac{c_{n}^2\eps\sqrt n}{4}.
\]
Note that
\[
|B_{n}(t)-B_{n}(t_{j})|
\leq
\sup_{t\in [L(0),L(1)]}\sup_{|t-y|\leq k_n^{-1}}|B_n(t)-B_n(y)|.
\]
Furthermore, for $(t,u)\in A_n'$ we have $t-2c_nu=L(x-2c_nz)\in J=[L(0),L(1)]$,
so that
\[
|B_{n}(t-2c_nu)-B_{n}(t_{j}-2c_nu)|\leq
\sup_{t\in J}\sup_{|t-y|\leq k_n^{-1}}|B_n(t)-B_n(y)|.
\]
Hence, from the triangle inequality it follows that
\[
\begin{split}
&
\sup_{t\in[t_{j-1},t_j]}\sup_{u\geq a\inf_tL'(t)}
\Big\{
B_{n}(t-2c_nu)-B_{n}(t_{j}-2c_nu)-B_{n}(t)+B_{n}(t_{j})
\Big\}\\
&\quad\leq
2\sup_{t\in J}\sup_{|t-y|\leq k_{n}^{-1}}
\left|B_n(t)-B_n(y)\right|.
\end{split}
\]
We conclude that
\[
\begin{split}
&
\prob\left(
\sup_{(x,z)\in A_n}
Y_n(x,z)>0\right)
\\
&\quad\leq
\prob\left(
2\sup_{t\in J}\sup_{|t-y|\leq k_{n}^{-1}}
\left|B_n(x)-B_n(y)\right|>\frac{c_{n}^2\eps\sqrt n}{4}
\right)\\
&\qquad+
\sum_{j=1}^{k_{n}}
\prob\left(\sup_{u\geq a\inf_{s}L'(s)}
\left\{
B_{n}(t_{j}-2c_nu)-B_{n}(t_{j})-\frac{c_{n}^2\eps \sqrt{n}u^2}{||L'||_{\infty}^2}
\right\}>0
\right).
\end{split}
\]
With Lemma~\ref{lem:increments WGen}, we have
\[
\begin{split}
&
\prob\Big(
2\sup_{t\in J}\sup_{|t-y|\leq k_{n}^{-1}}
\left|B_n(t)-B_n(y)\right|>\frac{c_{n}^2\eps\sqrt n}{4}
\Big)\\
&\quad\leq
K_{1}k_n\exp\left(-\frac{K_2\eps^2n}{16}c_{n}^4 k_{n}^\tau\right)
\leq
K_{1}\left(\frac{n}{c_{0}\log n}\right)^{1/(4-\tau)}n^{-K_2\eps^22^{-4-\tau}c_{0}},
\end{split}
\]
by definition of $c_{n}$ and $k_{n}$, since $ k_{n}\leq c_n^{-1}$ and $k_{n}\geq c_{n}^{-1}/2$ for sufficiently large~$n$.
Hence, there exist positive numbers $K_{1}$ and $K_{2}$ that depend only on $f$, $L$ and~$C_{0}$ such that
$$
\prob\left(
2\sup_{t\in J}\sup_{|t-y|\leq k_{n}^{-1}}
\left|B_n(t)-B_n(y)\right|>\frac{c_{n}^2\eps\sqrt n}{4}
\right)
\leq K_{1}n^{-K_{2}c_{0}},
$$
for all $n$, provided $c_{0}\geq C_{0}$ for some sufficiently large $C_{0}$.
Furthermore, with~{\Bnn} we have
\[
\begin{split}
\sum_{j=1}^{k_{n}}
&\prob\left(\sup_{u\geq a\inf_{t}L'(t)}
\left\{
B_{n}(t_{j}-2c_nu)-B_{n}(t_{j})-\frac{c_{n}^2\eps \sqrt{n}u^2}{||L'||_{\infty}^2}
\right\}>0
\right)\\
&\leq
\sum_{j=1}^{k_{n}}
\prob\left(\sup_{z\geq 2c_{n}a\inf_{t}L'(t)}
\left\{
B_{n}(t_{j}-z)-B_{n}(t_{j})-\frac{\eps \sqrt n z^2}{4||L'||_{\infty}^2}
\right\}>0
\right)\\
&\leq \sum_{j=1}^{k_{n}} K_{1}
\exp\left(
-\frac{K_{2}\eps^2n}{16||L'||_{\infty}^4}\left(2c_{n}a\inf_{t}L'(t)\right)^{4-\tau}
\right)
\\
&\leq K_{1}\left(\frac{n}{c_{0}\log n}\right)^{1/(4-\tau)}n^{-K_{2}\eps^2c_{0}(2a\inf_{t}L'(t))^{4-\tau}/(16||L'||_{\infty}^4)},
\end{split}
\]
by definition of $c_{n}$ and $k_{n}$. Renaming $K_{1}$ and $K_{2}$, the right hand term in the previous display is bounded from above by $K_{1}n^{-K_{2}c_{0}}$ for all $n$, provided $c_{0}\geq C_{0}$ for some sufficiently large $C_{0}$, where $K_{1}$ and $K_{2}$ depend only on $f$, $L$ and $C_{0}$. We conclude that there exist $K_{1},K_{2}$ such that
$$
\prob\left(
\sup_{(x,z)\in A_n}
Y_n(x,z)>0\right)
\leq
K_{1}n^{-K_{2}c_{0}},
$$
for all $n$, provided $c_{0}$ is sufficiently large.
Using~\eqref{KWP12Gen},
we conclude that
\begin{equation}\label{eq:KWP12bisGen}
\prob_{1,2}
\leq
K_{1}n^{-K_{2}c_{0}}
+
\prob\left(
\sup_{x\in[2c_n,1]}\sup_{z\in[1,a]}Y_n(x,z)>0\right).
\end{equation}
Next, we consider the case $z\in[1,a]$ and
establish an upper bound for
the probability on the right hand side of~\eqref{eq:KWP12bisGen}. Since $c_{n}$ tends to zero as $n\to\infty$ and~$f'$ is uniformly continuous on $[0,1]$,
we have
$$
F(x-2c_nz)-F(x)+2c_{n}\alpha_{x}
z=2c_{n}^2|f'(x)|z(1-z)+o(c_{n}^2),
$$
where $2c_{n}^2|f'(x)|z(1-z)\leq 0$ and  $o(c_{n}^2)$ is uniform in $z\in[1,a]$ and $x\in[2c_n,1]$.
Therefore,
\begin{equation}
\label{eqexpG2Gen}
F(x-2c_nz)-F(x)+2c_{n}\alpha_{x} z+t_{x}\leq -c_{n}^2|f'(x)|/8,
\end{equation}
for all $z\in[1,a]$ and $x\in[2c_{n},1]$ provided that $n$ is sufficiently large. With $M_{n}^B=F_{n}^B-F$,  it follows from \eqref{eq:KWP12bisGen} that
\[
\begin{split}
&\prob_{1,2}
\leq
K_1 n^{-K_{2}c_{0}}\\
&\quad+
\prob\left(
\sup_{x\in[2c_n,1]}\sup_{z\in[1,a]}
\left(
M_{n}^B(x-2c_nz)-M_{n}^B(x)
\right)>
\frac{c_{n}^2}8\inf_{t\in[0,1]}|f'(t)|\right).
\end{split}
\]
Repeating the same arguments as above yields
$\prob_{1,2}\leq K_1 n^{-K_{2}c_{0}}$,
for some positive $K_{1},K_{2}$ that depend only on $f$, $L$ and $C_{0}$, for all $n$, provided that $c_{0}\geq C_{0}$ for some  sufficiently large $C_{0}$.

We have already proved that $\prob_1\leq \prob_{1,1}+\prob_{1,2}$,
where $\prob_{1,1}$ satisfies~\eqref{eq:bound P11Gen},
so from~\eqref{eq:KWtildeGen}, we derive that for some positive $K_{1}$ and $K_{2}$ we have
$$
\prob\left(\tilde x_i>x  \text{ for some } x\in[2c_n,1]\right)\leq \prob_2+K_1 n^{-K_{2}c_{0}}.
$$
To deal with $\prob _{2}$, one can write
$\prob_2\leq \prob_{2,1}+\prob_{2,2}$, where
$$
\prob_{2,1}
=
\prob\left(\exists x\in[2c_n,1]:\ F_n^B(x-2c_{n})-F_n^B(x-c_n)>-c_{n}\alpha_x+t_x \right)
$$
and
$
\prob_{2,2}$ is the probability that there exist
$x\in[2c_n,1]$ and $z\in[x,(x+2c_n)\wedge1]$, such that
$$F_n^B(z)-F_n^B(x-2c_n)
>
\left(z-x+2c_{n}\right)\alpha_x-t_x,
$$
where we recall that $\alpha_x=f( x)+c_n|f'(x)|$ and $t_x=c_{n}^2 f'(x)/4$.
One can then conclude, using similar arguments as above, that
there exist positive numbers~$K_{1}$ and~$K_{2}$ such that
$\prob_{2}\leq K_1 n^{-K_{2}c_{0}}$,
for all $n$, provided that $c_{0}$ is sufficiently large, whence \eqref{PxiGen}.
This concludes the proof of  Lemma \ref{lem:espGunifGen}. \tqed

\paragraph{Proof of Theorem \ref{theorem:KWGenB}.}
Recall that we only need to prove the theorem for the case $[a,b]=[0,1]$. In the sequel, for all intervals $I\subset\R$, we denote by $\mbox{CM}_{I}$ the operator that maps a bounded function $h:I\to\R$ into the least concave majorant of $h$ on $I$.
First note that  for all $x\in[0,1]$, we have
\begin{equation}\notag
\hat F_{n,c_n}^{(B,x)}(x)-F_{n}^B(x)=(\mbox{CM}_{I_{n}(x)}T_{n}^{(B,x)})(0),
\end{equation}
where $I_n(x)=\left[- \left((c_{n}^{-1}x)\wedge2\right),(c_{n}^{-1}(1-x))\wedge2\right]$ and
$$T_{n}^{(B,x)}(\eta)=F_{n}^B(x+c_{n}\eta)-F_{n}^B(x)$$
for all $\eta\in I_{n}(x).$
With $M_{n}^B=F_{n}^B-F$, we can write
\begin{equation}\label{eq:TnBx}
T_{n}^{(B,x)}(\eta)
=
M_{n}^B(x+c_{n}\eta)-M_{n}^B(x)
+
F(x+c_{n}\eta)-F(x).
\end{equation}
Since $f'$ is bounded and $|\eta|\leq 2$, for all $\eta\in I_{n}(x)$, it follows from Taylor's expansion that
\begin{equation}
\label{eq:KWtaylorGen}
T_{n}^{(B,x)}(\eta)
=
M_{n}^B(x+c_{n}\eta)-M_{n}^B(x)+Y_{n}^{(B,x)}(\eta)+O(c_{n}^2),
\end{equation}
where $Y_{n}^{(B,x)}(\eta)=c_{n}\eta f(x)$, and
where the big $O$-term is uniform in $\eta\in I_n(x)$ and $x\in[0,1]$.
Because the process~$Y_{n}^{(B,x)}$ is linear, its least concave majorant on $I_n(x)$ is $Y_{n}^{(B,x)}$ itself.
Using that the supremum distance between the least concave majorants of processes is less than or equal
to the supremum distance between the processes themselves, we conclude from~\eqref{eq:KWtaylorGen}
that for all $x\in[0,1]$,
\[
\begin{split}
|\hat F_{n,c_n}^{(B,x)}(x)-F_n^B(x)|
&=
|(\mbox{CM}_{I_{n}(x)}T_{n}^{(B,x)})(0)|\\
&\leq
|Y_{n}^{(B,x)}(0)|
+
|(\mbox{CM}_{I_{n}(x)}T_{n}^{(B,x)})(0)-Y_{n}^{(B,x)}(0)|\\
&\leq
 \sup_{\eta\in I_n(x)}
|T_{n}^{(B,x)}(\eta)-Y_{n}^{(B,x)}(\eta)|\\
&\leq
\sup_{\eta\in I_n(x)}\left|M_{n}^B(x+c_{n}\eta)-M_{n}^B(x)\right|+
O(c_{n}^{2}).
\end{split}
\]
Hence,   for $A>0$ sufficiently large we have
\[\begin{split}
&\prob\left(\sup_{x\in[0,1]}
|\hat F_{n,c_n}^{(B,x)}(x)-F_n^B(x)|
>Ac_{n}^2\right)\\
&\quad\leq
\prob\left(
\sup_{x\in[0,1]}
\sup_{\eta\in I_n(x)}\left|M_{n}^B(x+c_{n}\eta)-M_{n}^B(x)\right|>Ac_{n}^2/2\right)\\
&\quad\leq
\prob\left(
\sup_{x\in [0,1]}
\sup_{\eta\in I_n(x)}\left|B_{n}\circ L(x+c_{n}\eta)-B_{n}\circ L(x)\right|>Ac_{n}^2\sqrt n/2\right)\\
&\quad\leq
\prob\left(
\sup_{x\in [L(0),L(1)]}\sup_{|x-y|\leq 2c_n||L'||_{\infty}}
\left|B_{n}(x)-B_{n}(y)\right|>Ac_{n}^2\sqrt n/2\right),
\end{split}
\]
since $|\eta|\leq 2$  for all $\eta\in I_{n}(x)$. We obtain from Lemma~\ref{lem:increments WGen} that for $A>0$ sufficiently large,
\begin{equation}
\label{eq:PboundLoc}
\begin{split}
&\prob\left(\sup_{x\in[0,1]}
|\hat F_{n,c_n}^{(B,x)}(x)-F_n^B(x)|
>Ac_{n}^2\right)\\
&\quad\leq
\frac{K_{1}}{2c_{n}||L'||_\infty}\exp\left(-K_{2}A^22^{-2-\tau}||L'||_{\infty}^{-\tau}nc_{n}^{4-\tau}\right)\\
&\quad\leq
\frac{K_{1}}{2||L'||_\infty}\left(\frac{n}{c_{0}\log n}\right)^{1/(4-\tau)}n^{-K_{2}A^22^{-2-\tau}||L'||_{\infty}^{-\tau}c_{0}}.
\end{split}
\end{equation}
The upper bound tends to zero as $n\to\infty$ provided that $A$ is sufficiently large, whence
$$\sup_{x\in[0,1]}
|\hat F_{n,c_n}^{(B,x)}(x)-F_n^B(x)|=O_{p}(c_n^2).$$
This completes the proof of Theorem \ref{theorem:KWGenB}. \tqed

\paragraph{Proof of Theorem~\ref{theorem:KWGenE}.}
Write
\[
\hat F_n- F_n
=
(\hat F_n-\hat F_n^B)
+
( F_n^B- F_n)
+
(\hat F_n^B-F_n^B).
\]
Since the supremum distance between least concave majorant processes is less than or equal to the supremum distance between the processes themselves,
the triangle inequality yields
\begin{equation}
\label{eq:BtoEGen}
\sup_{x\in[a,b]}|\hat F_n(x)- F_n(x)|
\leq
2\sup_{x\in[a,b]}|F_n(x)- F_n^B(x)|
+
\sup_{x\in[a,b]}|\hat F_n^B(x)-F_n^B(x)|.
\end{equation}
Theorem~\ref{theorem:KWGenE} now follows from assumption~\eqref{eq:embedGen}
and Theorem~\ref{theorem:KWGenB}.
\tqed

\paragraph{Proof of Theorem~\ref{theorem:KWSBGenEsp}.}
Recall that we only have to prove the theorem for the case $[a,b]=[0,1]$. In the sequel, we use the same notation as in Lemma \ref{lem:espGunifGen} and consider an arbitrary $r\geq 1$.
It follows from Fubini's Theorem that for any $A_{0}>0$ we have
\[
\begin{split}
&\E\left[\sup_{x\in[0,1]}|\hat F_n^B(x)-F_{n}^{B}(x)|^r\right]\\
&\quad=
\int_{0}^\infty \prob\left(\sup_{x\in[0,1]}|\hat F_n^B(x)-F_{n}^{B}(x)|^r>u\right)\,\mathrm{d} u\\
&\quad\leq
\left(A_{0}c_{n}^{2}\right)^r
+
\int_{A_{0}c_{n}^{2}}^\infty \prob\left(\sup_{x\in[0,1]}|\hat F_n^B(x)-F_{n}^{B}(x)|>v\right)rv^{r-1}\,\mathrm{d}v,
\end{split}
\]
where we used the fact that a probability is smaller than or equal to one, and we performed
a change of variable $v=u^{1/r}$.
From the triangle inequality, it follows that
\begin{equation}
\label{eq:In1In2}
\E\left[\sup_{x\in[0,1]}|\hat F_n^B(x)-F_{n}^{B}(x)|^r\right]
\leq
\left(A_{0}c_{n}^{2}\right)^r+ I_{n1}+I_{n2},
\end{equation}
where
$$
I_{n1}=\int_{0}^\infty \prob\left(
\sup_{x\in[0,1]}|\hat
F_{n,c_n}^{(B,x)}(x)-\hat F_n^B(x)|>\frac{v}{2}\right)rv^{r-1}\,\mathrm{d}v
$$
and
$$
I_{n2}=\int_{A_{0}c_{n}^{2}}^\infty\prob\left(
\sup_{x\in[0,1]}|\hat
F_{n,c_n}^{(B,x)}(x)-F_n^B(x)|>\frac{v}2\right)rv^{r-1}\,\mathrm{d}v.
$$

First, consider $I_{n1}$. It follows from Lemma  \ref{lem:espGunifGen}  that there exist positive $K_{1}$ and $K_{2}$ such that for all $v>0$,
\begin{equation}
\label{eq:majlemme1Gen}
\prob\left(
 \sup_{x\in[0,1]}|\hat
F_{n,c_n}^{(B,x)}(x)-\hat F_n^B(x)|>\frac{v}{2}\right)\leq K_{1} n^{-K_{2}c_{0}}.
\end{equation}
Moreover, similar to~\eqref{eq:BtoEGen}, the triangle inequality yields
\[
\begin{split}
\sup_{x\in[0,1]}|\hat
F_{n,c_n}^{(B,x)}(x)-\hat F_n^B(x)|
&\leq
\sup_{x\in[0,1]}|\hat
F_{n,c_n}^{(B,x)}(x)-F(x)|+\sup_{x\in[0,1]}|
F(x)-\hat F_n^B(x)|\\
&\leq
2 \sup_{x\in[0,1]}|
F_{n}^B(x)-F(x)|.
\end{split}
\]
By definition, $F_{n}^B-F=n^{-1/2}B_{n}\circ L$, so together with~\eqref{eq:BnSup} we derive that for all $v>0$,
\begin{equation}
\label{eq:maj2lemme1}
\begin{split}
\prob\left(
\sup_{x\in[0,1]}|\hat
F_{n,c_n}^{(B,x)}(x)-\hat F_n^B(x)|>\frac{v}{2}\right)
&\leq
\prob\left(
\sup_{x\in[L(0),L(1)]}|B_{n}(x)|>\frac{v\sqrt n}{4}\right)\\
&\leq
K_{1}\exp(-K_{2}nv^2/16).
\end{split}
\end{equation}
Note that without loss of generality, possibly enlarging $K_{1}$ and diminishing $K_{2}$,
we can choose $K_{1}$ and $K_{2}$ to be the same in~\eqref{eq:majlemme1Gen} and~\eqref{eq:maj2lemme1}.
Using the bound~\eqref{eq:majlemme1Gen} for $v\leq n$ and the bound~\eqref{eq:maj2lemme1} for $v>n$, we obtain
\[\begin{split}
&I_{n1}\leq \int_{0}^{n}K_{1} n^{-K_{2}c_{0}}rv^{r-1}\,\mathrm{d}v+\int_{n}^\infty K_{1}\exp(-K_{2}nv^2/16)rv^{r-1}\,\mathrm{d}v.
\end{split}\]
Consider a number $q$ with
\begin{equation}
\label{eq:q}
r-3q<-2r/(4-\tau)
\end{equation}
and let $K_{3}=q^q\exp(-q)$, so that $x^q\exp(-x)\leq K_{3}$ for all $x\in[0,\infty)$. Then,  we arrive at
\begin{equation}
\label{eq:In1}
\begin{split}
I_{n1}
&\leq
K_{1}n^{-K_{2}c_{0}+r}+K_{1}K_{3}\int_{n}^\infty (K_{2}nv^2/16)^{-q}rv^{r-1}\,\mathrm{d}v \\
&\leq
K_{1}n^{-K_{2}c_{0}+r}+K_{1}K_{3}\left(K_{2}n/16\right)^{-q}\frac{rn^{r-2q}}{2q-r},
\end{split}
\end{equation}
since $r-2q<0$ for all $q$ that satisfy \eqref{eq:q}.
Choose $c_{0}$ sufficiently large so that
$-K_{2}c_{0}+r<-2r/(4-\tau)$.
Then from \eqref{eq:q} and \eqref{eq:In1}, we conclude that
\begin{equation}\label{eq:boundIn1}
I_{n1}=o\left(\frac{\log n} n\right)^{2r/(4-\tau)}.
\end{equation}

Next, consider $I_{n2}$. Using a change of variable, we have
$$
I_{n2}
=
\left(2c_n^2\right)^r
\int_{A_{0}/2}^\infty\prob\left(
\sup_{x\in[0,1]}|\hat
F_{n,c_n}^{(B,x)}(x)-F_n^B(x)|>c_{n}^2v\right)rv^{r-1}\,\mathrm{d}v.
$$
Then we derive from \eqref{eq:PboundLoc}, that for sufficiently large $A_{0}$,
$$
I_{n2}
\leq
\left(2c_n^2\right)^r
\frac{K_{1}}{2||L'||_\infty}
\int_{A_{0}/2}^\infty
\left(\frac{n}{c_{0}\log n}\right)^{1/(4-\tau)}
n^{-K_{2}v^22^{-2-\tau}||L'||_{\infty}^{-\tau}c_{0}}rv^{r-1}\,\mathrm{d}v.
$$
Let $q>0$ with \eqref{eq:q} and let $K_{3}=q^q\exp(-q)$.
Then, similar to~\eqref{eq:In1}, the integral in
the previous display is bounded from above by
$$
K_{3}\int_{A_{0}/2}^\infty (c_{0}\log n)^{-1/(4-\tau)}
\left[\left(K_{2}v^22^{-2-\tau}||L'||_{\infty}^{-\tau}c_{0}-\frac{1}{4-\tau}\right)\log n\right]^{-q} rv^{r-1}\,\mathrm{d}v.
$$
Choosing $A_{0}$ sufficiently large,
so that
$$
K_{2}(A_{0}/2)^22^{-2-\tau}||L'||_{\infty}^{-\tau}c_{0}>\frac{2}{4-\tau},
$$
this is bounded from above by
$$K_{3}\int_{A_{0}/2}^\infty (c_{0}\log n)^{-1/(4-\tau)}
\left[
K_{2}v^22^{-3-\tau}||L'||_{\infty}^{-\tau}c_{0}\log n
\right]^{-q} rv^{r-1}\,\mathrm{d}v.
$$
Hence, $I_{n2}$ is bounded by
$$
\left(2c_n^2\right)^r
\frac{K_{1}K_{3}}{2||L'||_\infty}(c_{0}\log n)^{-1/(4-\tau)}\left[K_{2}2^{-3-\tau}||L'||_{\infty}^{-\tau}c_{0}\log n\right]^{-q} \frac{r(A_{0}/2)^{r-2q}}{2q-r},
$$
since $r-2q<0$ for all $q$ with \eqref{eq:q}. We conclude that for $A_{0}$ sufficiently large,
\begin{equation}\label{eq:boundIn2}
I_{n2}=o\left(\frac{\log n} n\right)^{2r/(4-\tau)}
\end{equation}
by definition of $c_{n}$. Combining \eqref{eq:In1In2}, \eqref{eq:boundIn1} and  \eqref{eq:boundIn2} completes the proof of Theorem~\ref{theorem:KWSBGenEsp}.
\tqed

\paragraph{Proof of Theorem~\ref{theorem:KWEGenEsp}.}

By convexity, we have $(a+b)^r\leq 2^{r-1}(a^r+b^r)$ for all positive numbers $a,b$ and therefore, \eqref{eq:BtoEGen} yields
\[
\begin{split}
&\E\left[\sup_{x\in[a,b]}|\hat F_n(x)-F_{n}(x)|^r\right]\\
&\quad\quad\leq2^{2r-1}\E\left[\sup_{x\in[a,b]}|F_n^B(x)-F_{n}(x)|^r\right]+2^{r-1}\E\left[\sup_{x\in[a,b]}|\hat F_n^B(x)-F_{n}^{B}(x)|^r\right].
\end{split}
\]
Theorem~\ref{theorem:KWEGenEsp}  then follows from \eqref{eq:moments} combined with Theorem \ref{theorem:KWSBGenEsp}. \tqed

\subsection{Proofs for Subsection \ref{sec:local}}

\paragraph{Proof of Theorem \ref{theorem:localrate}.}
In the sequel, we use the same notation as in Lemma~\ref{lem:espGunifGen}.
We first prove that
\begin{equation}\label{eq:theolocalB}
\sup_{|x-x_{0}|\leq\eps_{n}}|\hat F_{n}^B(x)-F_{n}^B(x)|=O_{p}\left((\eps_{n}^{\tau/2}n^{-1/2})\wedge\left(\frac{\log n}{n}\right)^{2/(4-\tau)}\right).
\end{equation}
It follows from Lemma \ref{lem:espGunifGen} that with $c_{0}$ sufficiently large,
$$
\sup_{|x-x_{0}|\leq\eps_{n}}|\hat F_n^B(x)-\hat F_{n,c_n}^{(B,x)}(x)|=O_{p}(\eps_{n}^{\tau/2}n^{-1/2}),
$$
so the triangular inequality yields
\begin{equation}\label{eq:triangle}
\begin{split}
&
\sup_{|x-x_{0}|\leq\eps_{n}}|\hat F_{n}^B(x)-F_{n}^B(x)| \\
&\leq
\sup_{|x-x_{0}|\leq\eps_{n}}|\hat F_{n,c_n}^{(B,x)}(x)-F_{n}^B(x_{0})-F(x)+F(x_{0})|\\
&\qquad+
\sup_{|x-x_{0}|\leq\eps_{n}}|F_n^B(x)- F_{n}^B(x_{0})-F(x)+F(x_{0})|
+
O_{p}(\eps_{n}^{\tau/2}n^{-1/2}).
\end{split}
\end{equation}
By definition, with $n$ sufficiently large and $x\in[x_{0}-\eps_{n},x_{0}+\eps_{n}]$, $\hat F_{n,c_n}^{(B,x)}$ is the least concave majorant of the restriction of $F_{n}^B$ to $I_{n}(x)=[x-2c_{n},x+2c_{n}]\cap[0,1]$ so by Marshall's lemma,
$$\sup_{\eta\in I_{n}(x)}|\hat F_{n,c_n}^{(B,x)}(\eta)-h(\eta)|\leq \sup_{\eta\in I_{n}(x)}|F_{n}^{B}(\eta)-h(\eta)|$$
for all concave functions $h:I_{n}(x)\to\R$.
The function $\eta\mapsto F_{n}^B(x_{0})+F(\eta)-F(x_{0})$ is concave on its domain, so Marshall's lemma ensures
that for all $x\in[x_{0}-\eps_{n},x_{0}+\eps_{n}]$,
$$|\hat F_{n,c_n}^{(B,x)}(x)-F_{n}^B(x_{0})-F(x)+F(x_{0})|\leq \sup_{\eta\in I_{n}(x)}|F_{n}^{B}(\eta)-F_{n}^B(x_{0})-F(\eta)+F(x_{0})|.$$
Setting $\delta_n=2c_{n}+\eps_{n}$, we conclude from \eqref{eq:triangle} that with $M_{n}^B=F_{n}^B-F$,
\begin{equation}
\label{eq:??}
\begin{split}
&
\sup_{|x-x_{0}|\leq\eps_{n}}|\hat F_{n}^B(x)-F_{n}^B(x)|\\
&\qquad\leq
2\sup_{|x-x_{0}|\leq \delta_{n}}|M_n^B(x)- M_{n}^B(x_{0})|+O_{p}(\eps_{n}^{\tau/2}n^{-1/2}).
\end{split}
\end{equation}
Consider the first term on the right hand side. With $y_{0}=L(x_{0})$  write
\begin{equation}
\label{????}
\begin{split}
n^{1/2}\sup_{|x-x_{0}|\leq\delta_{n}}|M_{n}^B(x)- M_{n}^B(x_{0})|
&=
\sup_{|x-x_{0}|\leq\delta_{n}}|B_{n}(L(x)) -B_{n}(L(x_{0}))|\\
&\leq
\sup_{|y-y_{0}|\leq\delta_{n}\|L'\|_{\infty}}|B_{n}(y) -B_{n}(y_{0})|,
\end{split}
\end{equation}
using that the derivative $L'$ is bounded. It follows from {\Bn} that for all $A>0$,
$$
\prob\left(\sup_{|y-y_{0}|\leq\delta_{n}\|L'\|_{\infty}}|B_{n}(y) -B_{n}(y_{0})|>A\delta_{n}^{\tau/2}\right)\leq K_{1} \exp(-K_{2} A^2\|L'\|_{\infty}^{-\tau}),
$$
which tends to zero as $A\to\infty$.
Combining this with the assumption that $\delta_{n}=2c_{n}+\eps_{n}\leq (2c_{0}^{1/(4-\tau)}+1)\eps_{n}$ yields
\begin{equation}
\label{eq:boundB}
\sup_{|y-y_{0}|\leq\delta_{n}\|L'\|_{\infty}}|B_{n}(y) -B_{n}(y_{0})|
=
O_{p}(\delta_{n}^{\tau/2})
=
O_{p}(\eps_{n}^{\tau/2}).
\end{equation}
Combining this with \eqref{????} and \eqref{eq:??} then yields
$$\sup_{|x-x_{0}|\leq\eps_{n}}|\hat F_{n}^B(x)-F_{n}^B(x)|=O_{p}(\eps_{n}^{\tau/2}n^{-1/2}).$$
Equation \eqref{eq:theolocalB} now follows from the previous display combined with Theorem~\ref{theorem:KWGenB}.
Similar to~\eqref{eq:BtoEGen}, we obtain
\begin{equation}
\label{eq:triangleLoc}
\begin{split}
\sup_{|x-x_{0}|\leq\eps_{n}}|\hat F_n(x)- F_n(x)|
&\leq
2\sup_{x\in[a,b]}|F_n(x)- F_n^B(x)|\\
&\qquad+
\sup_{|x-x_{0}|\leq\eps_{n}}|\hat F_n^B(x)-F_n^B(x)|.
\end{split}
\end{equation}
Hence, the theorem follows from~\eqref{eq:embedGen} and~ \eqref{eq:theolocalB}.
\tqed

\paragraph{Proof of Theorem \ref{lem:pointwiserate}.}
Combining~\eqref{eq:embedGen}, with $\gamma_n=O(n^{-2/(4-\tau)})$, and~\eqref{eq:triangleLoc}, with $\eps_{n}$ replaced by 0,
yields
\[
|\hat F_n(x_{0})- F_n(x_{0})|\leq
O_{p}(n^{-2/(4-\tau)})+|\hat F_n^B(x_{0})-F_n^B(x_{0})|.
\]
Therefore, it suffices to show that
\begin{equation}
\label{eq:pointwiserateB}
\hat F_{n}^B(x_{0})-F_{n}^B(x_{0})=O_{p}(n^{-2/(4-\tau)}).
\end{equation}
The proof of this is along the lines of the proof of Lemma~\ref{lem:espGunifGen},
except that we now take
\begin{equation}\label{eq: c0local}
c_{n}=(c_{0}/n)^{1/(4-\tau)},
\end{equation}
for some positive number $c_{0}$.
Without loss of generality we assume that $[a,b]=[0,1]$.
Define
$$
\tilde x_{0i}=\inf
\left\{u\geq  (x_{0}-2c_n)\vee 0\text{, such that }\hat F_n^B(u)= \hat F_{n,c_n}^{(B,x_{0})}(u)\right\},
$$
with the convention that the infimum of an empty set is $(x_{0}+2c_n)\wedge1$, and
$$
\tilde x_{0s}=
\sup
\left\{u\leq  (x_{0}+2c_n)\wedge 1\text{, such that }\hat F_n^B(u)= \hat F_{n,c_n}^{(B,x_{0})}(u)\right\},
$$
with the convention that the supremum of an empty set is $(x_{0}-2c_n)\vee0$.
Arguing as in the proof of Lemma \ref{lem:espGunifGen}, we obtain that
$$
\prob(\hat F_n^B(x_{0})\neq \hat F_{n,c_n}^{(B,x_{0})}(x_{0}))\leq \prob(\tilde x_{0i}> x_{0})+\prob(\tilde x_{0s}< x_{0}).
$$
Consider $\prob(\tilde x_{0i}> x_{0})$.
Note that if $x_0\leq 2c_n$, then $ \hat F_{n,c_n}^{(B,x_{0})}(0)=F_n^B(0)=\hat F_n^B(0)$, so that $\widetilde{x}_{0i}=0$
and hence, $\prob(\tilde x_{0i}> x_{0})=0$.
Next, consider $\prob(\tilde x_{0i}> x_{0})$, for $x_0\in[2c_n,1]$.
Let $\alpha_x=f(x_{0})+c_n|f'(x_{0})|$ and $t_x=c_{n}^2 f'(x_{0})/4$. Similar to the proof of Lemma \ref{lem:espGunifGen}, we have
$
\prob\left(\tilde x_{0i}>x_0\right)\leq \prob_1+\prob_2,$
where
$$
\prob_1=
\prob\Big(
 \exists 0\leq y\leq x_{0}-2c_n:
 F_n^B(y)-F_n^B(x_{0}-c_n)>(y-x_{0}+c_n)\alpha_x
\Big)
$$
and
$$
\prob_2=
\prob\Big(
\exists z\in[x_{0},(x_0+2c_n)\wedge1]:
F_n^B(z)-F_n^B(x_{0}-c_n)>(z-x_{0}+c_n)\alpha_x\Big).
$$
Furthermore,  $\prob_1\leq \prob_{1,1}+\prob_{1,2}$ where, with $\eps=\inf_{t}|f'(t)|$ and $K_{1},K_{2}$ as in Assumption {\Bn}, we have for sufficiently large $n$
\[
\begin{split}
\prob_{1,1}
&=
\prob\left(F_n^B(x_{0})-F_n^B(x_{0}-c_n)>c_{n}
\alpha_x+t_x \right)\\
&\leq
\prob\left(M_{n}^B(x_{0})-M_{n}^B(x_{0}-c_n)>\frac{c_{n}^2}{8}\inf_{t\in[0,1]}|f'(t)|\right)\\
&\leq
\prob\left(\sup_{|L(x_0)-y|\leq c_n||L'||_{\infty}}(B_{n}(L(x_0))-B_{n}(y))>\frac{c_{n}^2\sqrt n}{8}\eps\right)\\
&\leq
K_{1}\exp\left(-\frac{K_{2}\eps^2}{64||L'||_{\infty}^\tau}nc_{n}^{4-\tau}\right)
\end{split}
\]
and
$$
\prob_{1,2}
=
\prob\left(\exists 0\leq  y\leq x_{0}-2c_n:
F_n^B(x_{0})-F_n^B(y)<(x_{0}-y)\alpha_x+t_x\right).
$$
Similar to the proof of Lemma \ref{lem:espGunifGen}, for $z\in[1,x_0/(2c_n)]$, define
$$
Y_n(x_{0},z)=
F_n^B(x_{0}-2c_nz)-F_n^B(x_{0})
+2c_{n}\alpha_xz+t_x,
$$
so that
$\prob_{1,2}=
\prob\left(\exists z\in[1,x_0/(2c_n)]:Y_n(x_{0},z)>0
\right)$.
With $a$ such that $a\eps>2\sup_{t}|f'(t)|$, using~\eqref{eq:zGeqa} in the case $z\in[a,x_0/(2c_n)]$ and~\eqref{eqexpG2Gen}
in the case $z\in[1,a]$ we arrive at
\[
\begin{split}
\prob_{1,2}
&\leq
\prob\left(
\sup_{z\in[a,x_0/(2c_n)]}
\left\{
M_{n}^B(x_{0}-2c_nz)-M_{n}^B(x_{0})-\eps c_{n}^2 z^2
\right\}>0
\right)\\
&\qquad+
\prob\left( \sup_{z\in[1,a]}\left(M_{n}^B(x_{0}-2c_nz)-M_{n}^B(x_{0})\right)>
\frac{c_{n}^2}8\inf_{t\in[0,1]}|f'(t)|\right)\\
&\leq
K_{1}\exp\left(-K_{2}nc_{n}^{4-\tau}\right)
=
K_{1}\exp\left(-K_{2}c_0\right),
\end{split}
\]
for some positive $K_{1},K_{2}$ that depend only on $f$ and $L$.  For the last inequality, we used both assumptions {\Bn} and {\Bnn}. We used the definition \eqref{eq: c0local} for the last equality.
We conclude that $\prob_{1}$ tends to zero as $c_{0}\to\infty$.
Similarly, one can obtain that $\prob_{2}$ and $\prob(\tilde x_{0s}< x_{0})$ converge to zero,
as $c_{0}\to\infty$.
Therefore, for all $\eps>0$ there exists $C_{0}>0$ such that
$$
\prob\left(\hat F_n^B(x_{0})\neq \hat F_{n,c_n}^{(B,x_{0})}(x_{0})\right)\leq \eps.
$$
provided that $c_{0}\geq C_{0}$.
Hence, to prove \eqref{eq:pointwiserateB}, it now suffices to prove that
\begin{equation}\label{eq:pointwiserateLoc}
\hat F_{n,c_n}^{(B,x_{0})}(x_{0})-F_{n}^B(x_{0})= O_{p}(n^{-2/(4-\tau)}),
\end{equation}
for arbitrary $c_{0}$.
To this end, first note that
\begin{equation}\label{eq:CM}
\hat F_{n,c_n}^{(B,x_{0})}(x_{0})-F_{n}^B(x_{0})=(\mbox{CM}_{I_{n}}T_{n}^{(B,x_{0})})(0),
\end{equation}
where  $I_{n}=[-((c_{n}^{-1}x_0)\wedge2,(c_{n}^{-1}(1-x_0))\wedge 2]$,
$T_{n}^{(B,x_{0})}$ is taken from~\eqref{eq:TnBx}, and for all intervals $I\subset\R$,
$\text{CM}_{I}$ denotes the operator that maps a bounded function $h:I\to\R$ into the least concave majorant of $h$ on $I$.
Using~\eqref{????} and \eqref{eq:boundB} with~$\delta_{n}$ replaced by $2c_{n}$, we conclude that
$$
T_{n}^{(B,x_{0})}(\eta)
=O_{p}(n^{-1/2}c_{n}^{\tau/2})+F(x_{0}+c_{n}\eta)-F(x_{0})$$
where the big $O_{p}$-term is uniform in $\eta\in I_{n}$. Next, by Taylor expansion we have
\[
T_{n}^{(B,x_{0})}(\eta)
=
O_{p}(n^{-1/2}c_{n}^{\tau/2})+O(c_{n}^2)+c_{n}\eta f(x_{0})
=
O_{p}(c_{n}^2)+c_{n}\eta f(x_{0}),
\]
by definition \eqref{eq: c0local} of $c_{n}$, where the big $O_{p}$-term is uniform in $\eta\in I_{n}$. The supremum distance between the least concave majorants of processes is less than or equal to the supremum distance between the processes themselves, so with $Y_{n}(\eta)=c_{n}\eta f(x_{0})$, we have
$$(\mbox{CM}_{I_{n}}T_{n}^{(B,x_{0})})(0)=O_{p}(c_{n}^2)+(\mbox{CM}_{I_{n}}Y_{n})(0). $$
Since the process $Y_{n}$ is linear, we have  $\mbox{CM}_{I_{n}}Y_{n}=Y_{n}$ and therefore,
$(\mbox{CM}_{I_{n}}Y_{n})(0)=Y_{n}(0)=0.$ We then conclude from \eqref{eq:CM} that
$$\hat F_{n,c_n}^{(B,x_{0})}(x_{0})-F_{n}^B(x_{0})=O_{p}(c_{n}^2).$$
This completes the proof of \eqref{eq:pointwiserateLoc} by definition of $c_{n}$.
\tqed

\subsection{Proofs for Subsection \ref{sec:monoGen}.}

\paragraph{Proof of Corollary~\ref{cor:GenMono}.}
According to Theorems \ref{theorem:KWGenE} and \ref{theorem:KWEGenEsp}, it suffices to prove that $B_{n}$ satisfies {\Bn}, {\Bnn} and \eqref{eq:BnSup} with $\tau=1$. Note that we can write
\begin{equation}
\label{eq:Bn and WnGen}
B_{n}(t)=W_{n}(t)-\xi_{n}t,
\qquad
\text{for }t\in[a,b],
\end{equation}
where $W_{n}$ is  Brownian motion and
$\xi_{n}\equiv0$, if $B_{n}$ is Brownian motion,
and $\xi_{n}\sim N(0,1)$ independent of $B_{n}$,
if $B_{n}$ is Brownian bridge.
Therefore, for all fixed $x$ and all $u\in(0,1]$ and $v>0$ we have
\[
\begin{split}
\prob\left(\sup_{|x-y|\leq u}|B_{n}(x)-B_{n}(y)|>v\right)
&\leq
\prob\left(\sup_{|x-y|\leq u}|W_{n}(x)-W_{n}(y)|>v/2\right)\\
&\qquad+
\prob\left(u|\xi_{n}|>v/2\right).
\end{split}
\]
Using change of origin and scaling properties of Brownian motion, since $u\leq 1$ we obtain
\[
\begin{split}
\prob\left(\sup_{|x-y|\leq u}|B_{n}(x)-B_{n}(y)|>v\right)
&\leq
\prob\left(\sqrt u \sup_{|x|\leq 1}|W_{n}(x)|>v/2\right)\\
&\qquad+
\prob\left(\sqrt u|\xi_{n}|>v/2\right).
\end{split}
\]
By Doob's inequality (see e.g.~Proposition 1.8 in~\cite{revuz-yor1991}), the first probability on the right hand side is bounded by $4\exp(-v^2/(8u))$. Moreover, the second probability on the right hand side is bounded by $\exp(-v^2/(8u))$, whence
$$\prob\left(\sup_{|x-y|\leq u}|B_{n}(x)-B_{n}(y)|>v\right)\leq 5\exp(-v^2/(8u)).$$
This proves that $B_{n}$ satisfies {\Bn} with $K_{1}=5$, $K_{2}=1/8$ and $\tau=1$.
We obtain~\eqref{eq:BnSup} from~{\Bn} for the special case of $y=0$ and $u=L(b)-L(a)$, using that $B_{n}(0)=0$
almost surely.

Now, consider \Bnn. For all $u\in(0,1]$, $v>0$, and all $x\in I$ we have
\[
\begin{split}
&
\prob\left(\sup_{z\geq u}\left\{B_{n}(x-z)-B_{n}(x)-v z^2\right\}>0\right)\\
&\leq
\prob\left(\sup_{z\geq u}\left\{W_{n}(x-z)-W_{n}(x)-\frac{v z^2}2\right\}>0\right)
+
\prob\left(\sup_{z\geq u}\left\{\xi_{n}z-\frac{v z^2}2\right\}>0\right),
\end{split}\]
where $W_{n}$ and $\xi_{n}$ are taken from \eqref{eq:Bn and WnGen}.
Changing origin in the Brownian motion yields
\[
\begin{split}
&
\prob\left(\sup_{z\geq u}\left\{B_{n}(x-z)-B_{n}(x)-v z^2\right\}>0\right)\\
&\quad\leq
\prob\left(\sup_{z\geq u}\left\{W_{n}(z)-\frac{v z^2}2\right\}>0\right)+\prob\left(\sup_{z\geq u}\left\{\xi_{n}z-\frac{v z^2}2\right\}>0\right).
\end{split}
\]
By (3.3) in \cite{durottocquet2003}, the first probability on the right hand side is bounded by $\exp(-v^2u^3/8)$.
Moreover,
\[
\prob\left(\sup_{z\geq u}\left\{\xi_{n}z-\frac{v z^2}2\right\}>0\right)
\leq
\prob\left(|\xi_{n}|>vu/2\right)
\leq
\mathrm{e}^{-(vu)^2/8}
\leq
\mathrm{e}^{-v^2u^3/8},
\]
since $u\leq 1$.
Therefore,
$$
\prob\left(\sup_{z\geq u}\left\{B_{n}(x-z)-B_{n}(x)-v z^2\right\}>0\right)\leq 2\exp(-v^2u^3/8),
$$
which proves that $B_{n}$ satisfies {\Bnn} with $K_{1}=2$, $K_{2}=1/8$ and $\tau=1$. This concludes the proof of Corollary~\ref{cor:GenMono}. \tqed

\paragraph{Proof of Corollary~\ref{cor:regression}.}
Similar to Theorem~5(ii) in~\cite{durot2007},
it can be proved that if $\E|\epsilon_{i}|^{3}<\infty$, then
\begin{equation}
\label{eq:ass.embed}
\prob\left\{
n^{2/3}
\sup_{t\in[a,b]}
\left|
F_{n}(t)-\E(F_{n}(t))-n^{-1/2}B_{n}\circ L_{n}(t)
\right|
>x
\right\}
\leq
Cx^{-3},
\end{equation}
for all $x>0$,  with $B_{n}$ a Brownian  motion and
$$L_{n}(t)=\frac{\E(\epsilon_{i})^{2}}{n}\sum_{i=1}^n\indicator(t_{i}\leq t).$$
This implies that
$$\sup_{t\in[a,b]}|F_{n}(t)-\E(F_{n}(t))-n^{-1/2}B_{n}\circ L_{n}(t)|=O_{p}(n^{-2/3}).
$$
With \eqref{eq:design} and ~{\mono} we have
$$\sup_{t\in[a,b]}|F(t)-\E(F_{n}(t))|=O(n^{-2/3})$$ and therefore,
$$\sup_{t\in[a,b]}|F_{n}(t)-F(t)-n^{-1/2}B_{n}\circ L_{n}(t)|=O_{p}(n^{-2/3}).
$$
Moreover, with \eqref{eq:design} and $L(t)=(t-a)\E(\epsilon_{i})^{2}/(b-a)$,
we have
$$
\sup_{t\in[a,b]}|L_{n}(t)- L(t)|\leq Mn^{-1/3},
$$
for some $M>0$.
Since the Brownian motion $B_{n}$ satisfies {\Bn} with $\tau=1$, with Lemma~\ref{lem:increments WGen} we have
\[
\begin{split}
&
P
\left(
\sup_{t\in[a,b]}|B_{n}\circ L_{n}(t)-B_{n}\circ L(t)|>v
\right)\\
&\quad\leq
P\left(
\sup_{x\in[L(a),L(b)]}
\sup_{|x-y|\leq Mn^{-1/3}}
|B_{n}(x)-B_{n}(y)|
> v
\right)\\
&\quad\leq
K_{1}n^{1/3}M^{-1}\exp
\left(
-K_2v^2n^{1/3}/M
\right).
\end{split}
\]
It follows that
$$
n^{-1/2}\sup_{t\in[a,b]}|B_{n}\circ L_{n}(t)-B_{n}\circ L(t)|=O_{p}(n^{-2/3}(\log n)^{1/2}).
$$
Hence, \eqref{eq:embedGen} holds with $\gamma_n=O(n^{-2/3}(\log n)^{1/2})$ and $L(t)=(t-a)\E(\epsilon_{i}^2)/(b-a)$,
and Corollary~\ref{cor:regression} follows from Corollary~\ref{cor:GenMono}.
\tqed

\paragraph{Proof of Corollary  \ref{cor:density}.}
 From the proof of Theorem~6 in~\cite{durot2007}, it can be seen that, due to the Hungarian embedding,~\eqref{eq:embedGen} holds with $L=F$, $B_{n}$ a Brownian bridge, and  $\gamma_n=O(n^{-2/3})$. Therefore, Corollary~\ref{cor:density} follows from Corollary~\ref{cor:GenMono}.
\tqed
\paragraph{\bf Proof of Corollary \ref{cor:censoring}.}
Similar to Theorem 3 in~\cite{durot2007}, it can be proved that~\eqref{eq:embedGen} holds   with $B_{n}$ a Brownian  motion,
$$
L(t)=\int_{0}^t\frac{f(u)}{(1-G(u))(1-H(u))}\,\mathrm{d}u,\quad t\in[0,1]
$$
and $\gamma_n=O(n^{-2/3})$. Therefore, Corollary~\ref{cor:censoring} follows from Corollary~\ref{cor:GenMono}.

\tqed

\subsection{Proof of Corollary~\ref{cor:GenPos}.}
Let $F_{n}^B=F+n^{-1/2}B_{n}\circ L$, with $L$ and $B_n$ taken from~\eqref{eq:embedGen}.
Define
\[
S_n^B(t)=S(t)+n^{-1/2}\widetilde{B}_n\circ L(t)=F_n^B(b)-F_n^B(t),
\]
where
$\widetilde{B}_n(t)=B_n(L(b))-B_n(L(t))$.
Furthermore, let
$$
\widetilde B_{n}^H(t)=\int_a^{t} \widetilde B_{n}\circ L(x)\,\mathrm{d}x,
\quad
t\in[a,b],
$$
and define
$$
H_{n}^B(t)=\int_{a}^t
S_{n}^B(x)\,\mathrm{d}x
=
H(t)+n^{-1/2}\widetilde{B}_{n}^H(t),
\quad
t\in[a,b].
$$
Assumption \eqref{eq:embedGen} with $\gamma_n=O(n^{-1}\log n)$ ensures that
\[
\begin{split}
\sup_{t\in[a,b]}|H_{n}(t)-H_{n}^B(t)|
&\leq
(b-a)\sup_{t\in[a,b]}|S_{n}(t)-S_{n}^B(t)|\\
&\leq
(b-a)
\left(
|F_n(b)-F_n^B(b)|
+
\sup_{t\in[a,b]}|F_{n}(t)-F_{n}^B(t)|
\right)\\
&\leq
2(b-a)\sup_{t\in[a,b]}|F_{n}(t)-F_{n}^B(t)|\\
&=
O_{p}\left(\frac{\log n}{n}\right),
\end{split}
\]
which means that Assumption \eqref{eq:embedGen} also holds with $\gamma_n=n^{-1}\log n$,
and $F_{n},F,B_{n}$ replaced by $H_{n},H,\widetilde{B}_{n}^H$, respectively,
and $L(t)=t$.
Clearly $H(a)=0$ and, since $f=-S'$ satisfies~{\positive}, $H$ is twice continuously differentiable with a decreasing first derivative $S$ that satisfies
$$
0<\inf_{t\in[a,b]}|S'(t)|\leq \sup_{t\in[a,b]}|S'(t)|<\infty.
$$
We prove below that both {\Bn} and {\Bnn} hold with $\tau=2$, $L(t)=t$, and $B_n$ replaced by $\widetilde{B}_{n}^H$.
Then, Corollary~\ref{cor:GenPos} immediately follows from  Theorem~\ref{theorem:KWGenE}.

By definition of $\widetilde{B}_{n}^H$, for all $u>0$, $v>0$ and $x\in [a,b]$ we have
\[
\begin{split}
\sup_{|x-y|\leq u}|\widetilde{B}_{n}^H(x)-\widetilde{B}_{n}^H(y)|
&\leq
\sup_{|x-y|\leq u}
\int_{x\wedge y}^{x\vee y}|B_{n}\circ L(b)-B_n\circ L(t)|\,\mathrm{d}t\\
&\leq
2u\sup_{t\in[L(a),L(b)]}|B_{n}(t)|.
\end{split}
\]
Hence,
\[
\begin{split}
&
\prob\left(\sup_{|x-y|\leq u}|\widetilde B_{n}^H(x)-\widetilde B_{n}^H(y)|>v\right)\\
&\quad\leq
\prob\left(u\sup_{t\in[L(a),L(b)]}|W_{n}(t)|>v/4\right)
+
\prob\left(u\|L\|_{\infty}|\xi_{n}|>v/4\right)
\end{split}
\]
where $W_{n}$ and $\xi_{n}$ are taken from \eqref{eq:Bn and WnGen}, and where
$\|L\|_{\infty}=\sup_{t\in[a,b]}|L(t)| <\infty$, by assumption.
Therefore,
\[
\begin{split}
&
\prob\left(\sup_{|x-y|\leq u}|\widetilde B_{n}^H(x)-\widetilde B_{n}^H(y)|>v\right)\\
&\quad\leq
\prob\left(u\sup_{t\in[L(a),L(b)]}|W_{n}(t)|>v/4\right)
+
\exp(-v^2u^{-2}\|L\|_{\infty}^{-2}/32).
\end{split}
\]
By symmetry and scaling properties of Brownian motion, the first probability on the right hand side satisfies
\begin{eqnarray*}
\prob\left(u\sup_{t\in [L(a),L(b)]}|W_{n}(t)|>v/4\right)
&\leq&
2\prob\left(u\sqrt {\|L\|_{\infty}}\sup_{t\in [0,1]}|W_{n}(t)|>v/4\right)\\
&\leq&
4\prob\left(u\sqrt {\|L\|_{\infty}}\sup_{t\in [0,1]}W_{n}(t)>v/4\right).\\
\end{eqnarray*}
By Doob's inequality (see e.g.~Proposition 1.8 in~\cite{revuz-yor1991}), this is bounded by $4\exp(-v^2u^{-2}\|L\|_{\infty}^{-1}/32)$,
whence
\[
\prob\left(\sup_{|x-y|\leq u}|\widetilde{B}_{n}^H(x)-\widetilde{B}_{n}^H(y)|>v\right)
\leq
5\exp\left(-v^2u^{-2}(\|L\|_{\infty}^{-2}\wedge 1)/32\right).
\]
This proves that  {\Bn} holds with $B_{n}$ replaced by $\widetilde{B}_{n}^H$, $\tau=2$, $K_{1}=5$ and $K_2=(\|L\|_{\infty}^{-2}\wedge 1)/32$.

Next, consider~{\Bnn}.
By definition of $\widetilde{B}_{n}^H$, for all  $u>0$, $v>0$ and $x\in [a,b]$ we have
\[
\begin{split}
&
\prob\left(\sup_{z\geq u}\left\{\widetilde{B}_{n}^H(x-z)-\widetilde{B}_{n}^H(x)-v z^2\right\}>0\right)\\
&\quad\leq
\prob\left(\sup_{z\geq u}
\left\{2z\sup_{t\in [L(a),L(b)]}|B_{n}(t)|-v z^2\right\}>0\right)\\
&\quad\leq
\prob\left(\sup_{t\in [L(a),L(b)]}|B_{n}(t)|>vu/2\right).
\end{split}
\]
Similar arguments as above yield that {\Bnn} holds true with $B_{n}$ replaced by~$\widetilde{B}_{n}^H$, $\tau=2$,  $K_{1}=5$ and $K_2=(\|L\|_{\infty}^{-2}\wedge 1)/32.$
This completes the proof of Corollary~\ref{cor:GenPos}.
\tqed

\subsection{Proofs for Section \ref{sec:application}.}

\paragraph{Proof of Lemma \ref{lem:properties f_ns}.}
Let $p_1,\dots,p_m$ be the jump sizes of  $\hat f_{n}$ at the points of jump $\tau_1<\dots<\tau_m\in(a,b]$.
Note that $\hat f_{n}(x)=\hat f_{n}(b)$, for all $x\in(\tau_{m},b]$, and that
for $i=1,2,\dots,m$,
$$
\hat f_{n}(t)=\hat f_{n}(\tau_{i})=\hat f_{n}(b)+\sum_{j=i}^m p_j,
$$
for all $t\in(\tau_{i-1},\tau_{i}]$,
where $\tau_0=a$.
Therefore, when we define $K_{h_n}(t)=h_{n}^{-1}K(t/h_{n})$, for $t\in\R$, then for $t\in[a+h_{n},b-h_{n}]$, we can write
\[
\begin{split}
\hat f_{ns}(t)
&=
\frac{1}{h_n}\int_a^b K\left(\frac{t-x}{h_n}\right)\hat f_{n}(x)\,\md x.\\
&=\sum_{i=1}^m\left\{\hat f_{n}(b)+\sum_{j=i}^m p_j\right\}\int_{\tau_{i-1}}^{\tau_{i}}K_{h_n}(t-x)\,\md x
+\hat f_{n}(b)\int_{\tau_{m}}^bK_{h_n}(t-x)\,\md x,
\end{split}
\]
This means that for all $t\in[a+h_{n},b-h_{n}]$,
\[
\begin{split}
\hat f_{ns}(t)
&=
\sum_{j=1}^m p_j\int_a^{\tau_j}K_{h_n}(t-x)\,\md x+\hat f_{n}(b)\int_{a}^bK_{h_n}(t-x)\,\md x\\
&=
\sum_{j=1}^m p_j\int_{(t-\tau_j)/h_{n}}^{(t-a)/h_{n}}K(u)\,\md u +\hat f_{n}(b)\int_{(t-b)/h_{n}}^{(t-a)/h_{n}}K(u)\,\md u.
\end{split}
\]
Using that $K$ is supported on $[-1,1]$,
together with the fact that $(t-a)/h_{n}\geq 1$ and $(t-b)/h_{n}\leq -1$, for all $t\in[a+h_{n},b-h_{n}]$,
we obtain~\eqref{eq:implement}.
Because $K\geq0$, we conclude that $\hat f_{ns}$ is non-increasing on $[a+h_{n},b-h_{n}]$.
In particular,  we have $\widehat f_{ns}'(a+h_n)\leq 0$ and $\widehat f_{ns}'(b-h_n)\leq 0$, so
it immediately follows from definition~\eqref{eq:tildefKosorok} that $\widehat f_{ns}$ is also non-increasing
on the intervals $[a,a+h_{n}]$ and
$[b-h_{n},b]$. Since $\hat f_{ns}$ is continuous, we conclude that $\widehat f_{ns}$ is non-increasing on the whole interval $[a,b]$.
\tqed

\paragraph{Proof of Lemma \ref{lem:compar}.}
Denoting $\hat f_{ns}^{(0)}=\hat f_{ns}$ and $\hat f_{ns}^{(1)}=\hat f_{ns}'$,
for $l=0,1$ we have
\begin{equation}
\label{eq:smoothgren}
\begin{split}
&\sup_{t\in[a+h_{n},b-h_{n}]}\left|
\widehat{f}_{ns}^{(l)}(t)
-
\widetilde{f}_{n}^{(l)}(t)
\right|\\
&\quad\quad\quad=
\left|
\frac{1}{h_n^{1+l}}\int (\widehat{F}_{n}(t-uh_n)-F_{n}(t-uh_n))K^{(1+l)}(u)\,\md u
\right|\\
&\quad\quad\quad\leq
\frac{1}{h_n^{1+l}}\sup_{s\in[a,b]}|\widehat{F}_{n}(s)-F_{n}(s)|
\int |K^{(1+l)}|(u)|\,\md u\\
&\quad\quad\quad=
O_{p}\left(h_{n}^{-(1+l)}n^{-2/3}(\log n)^{2/3}\right),
\end{split}
\end{equation}
where we use Corollary~\ref{cor:GenMono} in the last equality.
On $[a,a+h_{n}]$ we have by definition~\eqref{eq:tildefKosorok},
\[\begin{split}
&\sup_{t\in[a,a+h_{n}]}\left|
\widehat{f}_{ns}(t)-\widetilde{f}_{n}(t)
\right|\\
&\quad\leq |\widehat f_{ns}(a+h_n)-\widetilde f_n(a+h_n)|+h_{n} |\widehat f_{ns}'(a+h_n)-\widetilde f_n'(a+h_n)|\\
&\quad\leq O_{p}\left(h_{n}^{-1}n^{-2/3}(\log n)^{2/3}\right),
\end{split}
\]
where we used~\eqref{eq:smoothgren} with $l=0,1$ in the last inequality.
Combining this with a similar argument on $[b-h_{n},b]$,
together with an application of~\eqref{eq:smoothgren} for $l=0$ on $[a+h_{n},b-h_{n}]$,
completes the proof of the lemma for $l=0$.
Similarly, for $l=1$,
\[
\begin{split}
\sup_{t\in[a,a+h_{n}]}
\left|
\widehat{f}_{ns}'(t)-\widetilde{f}_{n}'(t)
\right|
&=
\left|\widehat f_{ns}'(a+h_n)-\widetilde f_n'(a+h_n)\right|\\
&=
O_{p}\left(h_{n}^{-2}n^{-2/3}(\log n)^{2/3}\right).
\end{split}
\]
Using a similar argument on $[b-h_{n},b]$, together with applying~\eqref{eq:smoothgren} for $l=1$ on $[a+h_{n},b-h_{n}]$,
completes the proof for $l=1$.
\tqed

\paragraph{Acknowledgement.} 
The authors would like to thank the associate editor and two anonymous referees for their comments and suggestions,
which substantially improved the earlier version of the paper.
Special thanks goes to referee 1 for pointing out paper~\cite{wang2007kiefer}.

\bibliographystyle{acm}
\bibliography{KW}

\end{document}